\documentclass[11pt,reqno,a4paper]{amsart}
\usepackage{a4}

\title[HSL tori in $\CP^2$]{The classification of Hamiltonian stationary Lagrangian 
tori in $\CP^2$ by their spectral data.}
\author{Richard Hunter}
\author{Ian McIntosh}
\address{Department of Mathematics\\ University of York\\ 
York YO10 5DD, UK}
\email{ian.mcintosh@york.ac.uk}
\date{September 29, 2010}

\newcommand{\Z}{\mathbb{Z}}

\newcommand{\C}{\mathbb{C}}
\newcommand{\Ct}{\mathbb{C}^\times}
\newcommand{\Ci}{\mathbb{C}_\infty}
\newcommand{\R}{\mathbb{R}}

\renewcommand{\P}{\mathbb{P}}
\newcommand{\CP}{\mathbb{CP}}

\newcommand{\caA}{\mathcal{A}}
\newcommand{\caB}{\mathcal{B}}
\newcommand{\caC}{\mathcal{C}}
\newcommand{\caD}{\mathcal{D}}
\newcommand{\caE}{\mathcal{E}}
\newcommand{\caF}{\mathcal{F}}
\newcommand{\caG}{\mathcal{G}}
\newcommand{\caH}{\mathcal{H}}
\newcommand{\caI}{\mathcal{I}}

\newcommand{\caK}{\mathcal{K}}
\newcommand{\caL}{\mathcal{L}}
\newcommand{\caM}{\mathcal{M}}
\newcommand{\caO}{\mathcal{O}}

\newcommand{\caR}{\mathcal{R}}
\newcommand{\caS}{\mathcal{S}}
\newcommand{\caU}{\mathcal{U}}

\newcommand{\fA}{\mathfrak{A}}

\newcommand{\su}{\mathfrak{su}}
\newcommand{\gl}{\mathfrak{gl}}
\newcommand{\fb}{\mathfrak{b}}
\newcommand{\fg}{\mathfrak{g}}
\newcommand{\fh}{\mathfrak{h}}
\newcommand{\fk}{\mathfrak{k}}
\newcommand{\fm}{\mathfrak{m}}
\newcommand{\fo}{\mathfrak{o}}

\newcommand{\fz}{\mathfrak{z}}
\newcommand{\Lgt}{\Lambda^\tau\fg^\C}
\newcommand{\Lg}{\Lambda^\tau\fg}

\newcommand{\LCG}{\Lambda^\tau(C_\epsilon,G)}
\newcommand{\LCGsig}{\Lambda^\sigma(C_\epsilon,G)}
\newcommand{\LCGun}{\Lambda(C_{\epsilon^2},G)}
\newcommand{\Lcg}{\Lambda^\tau(C_\epsilon,\fg)}

\newcommand{\LIG}{\Lambda^\tau(I_\epsilon,G)}
\newcommand{\Lig}{\Lambda^\tau(I_\epsilon,\fg)}
\newcommand{\LIB}{\Lambda^\tau_B(I_\epsilon,G)}
\newcommand{\Lib}{\Lambda^\tau_\fb(I_\epsilon,\fg)}
\newcommand{\LEG}{\Lambda^\tau(E_\epsilon,G)}
\newcommand{\Leg}{\Lambda^\tau(E_\epsilon,\fg)}
\newcommand{\LHG}{\Lambda^\tau(\C^\times,G)}
\newcommand{\Lhg}{\Lambda^\tau(\C^\times,\fg)}
\newcommand{\Ltg}{\Lambda^\epsilon_{-2,\infty}}

\newcommand{\Aut}{\operatorname{Aut}}
\newcommand{\Lag}{\operatorname{Lag}}

\newcommand{\Hom}{\operatorname{Hom}}
\newcommand{\Pic}{\operatorname{Pic}}

\newcommand{\Jac}{\operatorname{Jac}}
\newcommand{\Prym}{\operatorname{Prym}}

\newcommand{\Spec}{\operatorname{Spec}}

\newcommand{\Res}{\operatorname{Res}}

\newcommand{\Sp}{\mathrm{Span}}
\newcommand{\Kah}{K\" ahler\ }
\newcommand{\pk}{polynomial Killing\ }
\newcommand{\HSL}{Hamiltonian stationary Lagrangian\ }
\newcommand{\Ad}{\operatorname{Ad}}
\newcommand{\Tr}{\operatorname{Tr}}
\newcommand{\ad}{\operatorname{ad}}

\newcommand{\e}{\epsilon}
\newcommand{\vol}{\mathrm{vol}}

\newtheorem{thm}{Theorem}[section]
\newtheorem{prop}[thm]{Proposition}
\newtheorem{lem}[thm]{Lemma}

\newtheorem{defn}[thm]{Definition}
\theoremstyle{remark}
\newtheorem{exam}{Example}
\newtheorem{rem}{Remark}

\numberwithin{equation}{section}

\begin{document}
\begin{abstract}
It is known that all weakly conformal \HSL 
immersions of tori in $\CP^2$ may be constructed by methods from integrable systems
theory. This article describes the precise details of a construction which leads to a form
of classification.  The immersion is encoded as spectral data in a
similar manner to the case of minimal Lagrangian tori in $\CP^2$, but the details 
require a careful treatment of both the ``dressing construction'' and the spectral data to deal
with a loop of flat connexions which is quadratic in the loop parameter. 
\end{abstract} 
\maketitle

\section{Introduction}

An immersion $f:M\to N$ of a compact Riemannian submanifold $M$ into 
a \Kah manifold $(N,\omega)$ is said to be Hamiltonian stationary Lagrangian\footnote{These 
are also known as H-minimal or Hamiltonian minimal.} (which we will abbreviate to HSL) if it is 
Lagrangian and its volume is 
critical for Hamiltonian variations. 
In \cite{Oh2} Oh derived the Euler-Lagrange equations in the following form.
Let $H$ be the mean curvature vector for $f$ and $\sigma_H$ the mean curvature 1-form
$f^*(H\rfloor\omega)$ on $M$. Oh showed that $M$ is Hamiltonian stationary if and only if 
$\sigma_H$ is co-closed. When $N$ is K\" ahler-Einstein
$\sigma_H$ is necessarily closed for any Lagrangian submanifold. Therefore a compact
Lagrangian immersion $f:M\to N$ is Hamiltonian stationary if and only if $\sigma_H$ is a harmonic
1-form. 

As well as being interesting in its own right, the HSL condition also arises in the study of 
minimal Lagrangian submanifolds of KE
manifolds (and not just because this is the case $H=0$). In their study of Lagrangian
surfaces in K\" ahler-Einstein 4-manifolds which minimise volume in their homotopy class, Schoen \& 
Wolfson \cite{SchW1} discovered that the minimiser may have isolated singularities whose 
tangent cone is only HSL. These are cones over HSL curves in $\CP^1$ and classification of the
latter allowed Schoen \& Wolfson to compute exactly which of these could give tangent cones for the
singular minimisers.

This provides a compelling reason for studying HSL submanifolds in $\CP^n$, since cones over these in
$\C^{n+1}$ provide the models for isolated conical singularities of Lagrangian volume minimisers in
higher dimensional KE manifolds. There are examples of HSL submanifolds of $\CP^n$ which are not 
minimal, mostly for the case $n=2$ (see, for example, \cite{CasHU,MaS,Mir1}) and all these examples are
highly symmetric. Because the study of minimal examples is much better developed throughout 
the rest of this article we will use ``HSL'' exclusively to mean non-minimal \HSL submanifolds.

The majority examples of HSL surfaces in $\CP^2$ are immersed tori (for an example of a HSL Klein 
bottle, see \cite{MirN}), 
and it is for this case that H\' elein \& Romon \cite{HelR3} (see also Ma \cite{Ma}) 
showed that the equations governing HSL tori lie in the realm of integrable systems. 
It was already known that for $H=0$ this approach gives a classification of minimal Lagrangian 
tori \cite{McI03} 
in terms of the associated spectral data (basically, an algebraic curve and a holomorphic
line bundle over it) and this led to an understanding of how
complicated the moduli space of minimal Lagrangian tori really is \cite{CarM,Has}. 

Our aim here is to provide the analogous classification by spectral data for HSL tori. To be
precise, let us define what we mean by HSL spectral data.
\begin{defn}\label{defn:spectral_data}
A triple $(X,\lambda,\caL)$ will be called spectral data for a (non-minimal) HSL torus in $\CP^2$
when the following conditions hold.
\begin{enumerate}
\item $X$ must be a complete connected algebraic curve possessing a holomorphic involution $\tau$ 
with only two fixed points and a real involution $\rho$, such that $\tau$ and $\rho$ commute.
In particular, the genus $g$ of $X$ is even.
\item $\lambda$ must be a rational function on $X$ of degree $3$, non-constant on every irreducible
component of $X$, satisfying 
\begin{equation}
\overline{\rho^*\lambda} = \bar\lambda^{-1},\quad \tau^*\lambda = -\lambda,
\end{equation}
and whose zeroes $P_1,P_2,P_3$ and poles $Q_1,Q_2,Q_3$ are distinct smooth points. $\tau$ must 
fix one zero $P_3$ (and hence one pole $Q_3$) of $\lambda$. $\rho$ must fix every point over the 
circle $|\lambda|=1$ and $\lambda$ must have no branch points on this circle. 
\item A double periodicity condition, which only depends upon $(X,\lambda)$, 
must be satisfied. The details of this are given in \S \ref{sec:Rec} below.
\item $\caL\simeq\caO_X(D)$ must be a line bundle of degree $g+2$ over $X$ and satisfy the linear 
equivalences 
\begin{equation}\label{eq:lequiv}
\rho_*D+D\sim R,\quad \tau_*D+D\sim R+P_3-Q_3,
\end{equation}
where $R$ is the ramification divisor of $\lambda$. Moreover, $D$ must contain no points on 
$|\lambda|=1$.
\end{enumerate}
\end{defn}
Now our principal result can be stated as follows.
\begin{thm}\label{thm:main}
There is an essentially bijective correspondence between:
\begin{enumerate}
\item weakly conformal, immersed, (non-minimal) HSL tori in $\CP^2$, and; 
\item triples of HSL spectral data $(X,\lambda,\caL)$. 
\end{enumerate}
\end{thm}
The word ``essentially'' means we work with the natural equivalence classes for such data.
Equations \eqref{eq:lequiv} can also be written in terms of the canonical class of $X$ (see
\eqref{eq:leqcan} below).

The real value of this theorem lies in the way the spectral data parameterises
the set of all HSL tori, and what this can tell us about families of HSL tori. The construction 
of explicit examples
is briefly discussed in \S\ref{sec:explicit}, where we show what must be done in the case where
$X$ is smooth, using the Riemann $\theta$-function. But such approaches are not very practical 
for curves of genus $4$ or more (and for genera $0$ and $2$, where only rational functions or
elliptic functions are required, the examples can be obtained faster
using the direct methods found in \cite{CasHU,HelR3,MaS,Mir1}).

The spectral data itself is best parameterised by the branch
divisor $B=\lambda_*R$ of $\lambda$, which determines $(X,\lambda)$ more or less uniquely as a
$3$-fold branched cover of the Riemann sphere $\Ci$. The symmetries $\rho,\tau$ will be consequences of
symmetries of $B$. In particular, the \emph{existence} of a HSL torus with spectral curve $X$
depends only on $(X,\lambda)$ and is independent of condition (d). In fact
(d) can be satisfied whenever $(X,\lambda)$ satisfy (a) and (b) (see \S\ref{sec:Rec} below). When 
equations \eqref{eq:lequiv} 
can be satisfied they can be satisfied by a real family of divisors of dimension $g/2$.
However, it is non-trivial matter to demonstrate 
that spectral data for non-minimal HSL tori exists for a given (even) spectral genus $g\geq 4$, 
although this will presumably follow from a perturbation of the argument used for minimal
Lagrangian tori in \cite{CarM}. Here we discuss our expectations for the dimension of the moduli 
space of HSL tori in \S\ref{sec:moduli}. 

The form of Theorem \ref{thm:main} looks very similar to the corresponding 
theorem for minimal Lagrangian tori \cite{McI03}, to which it reduces when the mean curvature is
zero (provided one interprets this limit carefully with respect to the double periodicity 
condition: see \S\ref{sec:moduli}). However, to achieve this required careful
adaptation of the tools used for that case, because the non-minimal HSL equations are encoded as 
the Maurer-Cartan
equations for a loop algebra valued $1$-form $\alpha_\zeta$ which has \emph{quadratic} dependence 
on the spectral parameter $\zeta$ (see equation \eqref{eq:MCform}). 
The success of the adaptation rested on two insights: (i) the natural vacuum solution for the 
relevant dressing action does \emph{not} correspond to a HSL immersion, but only gives a
constant map into $\CP^2$, (ii) the loop algebra can be ``partially untwisted'' (see
\S\ref{sec:untwist}) to enable us to view $\alpha_\zeta$ as linear in $\zeta^2$. 

To avoid assumptions about the
smoothness or irreducibility of the spectral curve we construct the spectral data from the
algebra of \pk fields (i.e, the spectral curve $X$ is a
ringed space and it carries a ``line bundle'' $\caL$ as an $\caO_X$-module), following
\cite{McI95,McI03}. In particular, our method of proof does not rule out the existence of
singular spectral curves. Indeed, this is a subtle point which is not well addressed in the
literature, since many spectral curve constructions force smoothness (by taking the normalisation) 
or assume it without comment. We prefer to allow for the possibility of singular spectral curves
since this simplifies, amongst other things, the periodicity conditions.

We expect the general method of our approach will apply to other examples where the Maurer-Cartan 
form is not linear in the spectral parameter (such as those in \cite{BurK,Khe,Ter}). In particular,
the appendix \ref{app:dressing} adapts the dressing construction of \cite{BurP2} to deal with the
general situation of a quadratic dependence on the spectral parameter.
By contrast, Mironov \cite{Mir2} gives an approach, specific to HSL, which uses a 
different form for the equations.

\smallskip\noindent
\textit{Acknowledgments.} This article includes and extends the work of the first author
which appeared in his dissertation \cite{Hun}. These results were first announced 
at the Symposium of the differential geometry of submanifolds at Valenciennes, July 2007
\cite{McIval}.  The approach taken here in defining the Lagrangian angle was suggested to us 
by Fran Burstall. We would also like to thank Fran Burstall and
David Calderbank for pointing out the inadequacies of the original attempt to 
describe the dressing arguments contained in appendix \ref{app:dressing}.

\smallskip\noindent
\textit{Notational conventions.} 
Throughout this article we denote $\C\setminus\{0\}$ by $\Ct$. We will use $e_1,\ldots,e_n$ to 
denote the standard (oriented orthonormal) basis
of $\R^n$ and $\C^n$. For any complex matrix $A$ its Hermitian transpose will be denoted by
$A^\dagger$. When $A$ is invertible we use $A^*$ to denote $(A^t)^{-1}$.
For any non-zero 
vector $v\in\C^{n+1}$ we denote by $[v]\in\CP^n$ the line generated by $v$. For a real Lie algebra
$\fg$ with complexification $\fg^\C$ conjugation with respect to that real form will be denoted by
$\bar\xi$ for $\xi\in\fg$. For any $1$-form $\beta\in\Omega^1_M$ its type
decomposition will be denoted by $\beta'+\beta''$, where $\beta'\in\Omega_M^{(1,0)}$ and
$\beta''\in\Omega_M^{(0,1)}$. For the most part the surfaces we deal with are tori, which we will
write as $\C/\Gamma$ (for a lattice $\Gamma$) when we want to emphasize the complex structure,
or as $\R^2/\Gamma$ when we want to emphasize the real structure. 

\section{The HSL condition via a loop of flat connexions.}

We begin by clarifying the notion of the \emph{Lagrangian angle} for Lagrangian submanifolds
of a K\" ahler-Einstein manifold. Let $(N,\omega)$ be a KE $2n$-manifold with K\" ahler form $\omega$
and volume form $\vol_N$. Let $\Lag(N)\to N$ be its bundle of oriented Lagrangian $n$-planes in $TN$. 
To any oriented Lagrangian immersion $f:M\to N$ we can assign the Gauss map 
\[
\gamma:M\to f^{-1}\Lag(N);\quad \gamma(p) = df(T_pM),
\]
and there is a unique smooth section $\Omega:M\to f^{-1}K_N$ for which at each point
$\Omega|\gamma(p)=\vol_N$: we will say $\Omega$ calibrates $\gamma(p)$. In that case
$\Omega$ must be a section of the unit circle subbundle $\caS$ of $f^{-1}K_N$. 
\begin{defn}
The unique section $\Omega:M\to \caS\subset f^{-1}K_N$ of the unit circle
bundle $\caS$ which calibrates the Gauss map $\gamma$ of $f$ at each point will be called
the \emph{Lagrangian angle} of $f$. 
\end{defn}
When $N$ is K\" ahler-Einstein $\caS$ is flat, since $K_N$
has curvature given, up to scale, by $\omega$. We define the Maslov form by
\begin{equation}\label{eq:Maslov}
\mu = -\frac{1}{i\pi}\Omega^*\alpha
\end{equation}
where $\alpha$ is the connexion 1-form for the induced connexion on 
the principal circle bundle $S$. Thus $\mu$ is closed, since $\caS$ is flat, and provides
the Maslov class $[\mu]\in H^1(M,\R)$.  A calculation originally due to Dazord \cite{Daz}
asserts that the mean curvature 1-form $\sigma_H = f^*(H\rfloor\omega)$ satisfies
\[
\sigma_H = \frac{\pi}{2}\mu.
\]
Thus $f$ is HSL if and only if it has co-closed, and therefore harmonic, Maslov form.

We focus now on the case where $N=\CP^2$.  We consider $\CP^2$ as the homogeneous space $G/K$
for the projective unitary group $G=PU(3)$ (i.e., $SU(3)$ modulo its centre), where $K$ is the
stabiliser of the line $[e_3]$. To make definitions easier, elements of $PU(3)$ will
tend to be thought of as represented by elements of
$SU(3)$.  Thus $K$ is the fixed point subgroup for the involution
$\sigma\in\Aut(G)$ defined by
\begin{equation}\label{eq:sigma}
\sigma(g) = \begin{pmatrix} -I_2 & 0 \\ 0 & 1\end{pmatrix} g \begin{pmatrix} -I_2 & 0 \\ 0 &
1\end{pmatrix},
\end{equation}
where $I_2$ is the $2\times 2 $ identity matrix. The corresponding reductive decomposition of
$\fg\simeq\su(3)$ will be written $\fg=\fk +\fm$. Recall that the tangent space
$T^\C(G/K)$ can be identified with the quotient space $G\times_K\fm^\C$, where $K$ 
acts adjointly on $\fm^\C$. 

Now define an order 4 outer automorphism $\tau\in\Aut(G)$ by
\begin{equation}\label{eq:tau}
\tau(g) = Sg^*S^{-1},\quad S = \begin{pmatrix} J & 0 \\ 0 & 1\end{pmatrix},\quad
J = \begin{pmatrix} 0 & 1\\ -1 & 0\end{pmatrix}.
\end{equation}
Since $\tau$ preserves the centre of $SU(3)$ it is well-defined on $PU(3)$. We will denote the fixed
point subgroup of $\tau$ by $G_0$: it is easy to check that $G_0\simeq SU(2)$. Of course, $G/G_0$ is a
4-symmetric space. 
It can be shown (see \cite{BurK,McIval}) that it is isomorphic to the unit subbundle 
of the canonical bundle $K_{\CP^2}$.

We are now in a position to introduce a preferred family of frames for Lagrangian immersions $f:M\to
G/K$. As we have explained above, any such immersion possesses a  Lagrangian angle $\Omega:M\to G/G_0$, 
thought of as a section of $f^{-1}K_{\CP^2}$. 
\begin{defn}
Let $f:M\to G/K$ be a Lagrangian immersion and $\Omega:M\to G/G_0$ its Lagrangian angle. A frame 
$F:\tilde M\to G$ will be called a \emph{Lagrangian frame} for $f$ if it also frames $\Omega$.
\end{defn}
Of course, Lagrangian frames are not unique: they may be right multiplied by any element of the gauge
group $\caG_0=C^\infty(\tilde M,G_0)$.

Given a Lagrangian frame $F$ its left Maurer-Cartan form $\alpha=F^{-1}dF$ admits a decomposition
$\alpha_\fk+\alpha_\fm$ into summands taking values in the eigenspaces of $\sigma$. 
Recall that
$\Ad F\cdot\alpha_\fm$ represents the tangent map $df$. Now let $\fg_j\subset \fg^\C$ denote the
$i^j$-eigenspace of $\tau$ (as the induced automorphism on $\fg^\C$) and notice that since 
$\sigma=\tau^2$ we must have $\fk^\C = \fg_0+\fg_2$, and
$\fm^\C = \fg_{-1}+\fg_1$. To be explicit, these subspaces $\fg_j\subset\mathfrak{sl}_3(\C)$ are
\begin{eqnarray}\label{eq:g_j}
& \fg_0 = \{\begin{pmatrix} Q & 0\\ 0&0\end{pmatrix}: Q\in\mathfrak{sl}_2(\C)\}, 
&\fg_1 = \{\begin{pmatrix} O_2& u\\-i(Ju)^t & 0 \end{pmatrix}: u\in\C^2\},\notag \\
&\fg_{-1} = \bar\fg_1,& \fg_2 = \{\begin{pmatrix} aI_2 & 0 \\ 0 & -2a\end{pmatrix}:a\in\C\}.
\end{eqnarray}
Here we are using $O_2$ to denote the $2\times 2$ zero matrix.
Note in particular that $[\fg_0,\fg_2]$=0, hence the adjoint action of $G_0$ on $\fg_2$ is trivial.

Now let $\alpha_j$ denote the projection of $\alpha$ onto $\fg_j$ for $j=-1,0,1$, then
$\alpha_\fm = \alpha_{-1}+\alpha_1$ and the condition that $f$ is conformal Lagrangian is that
$\alpha^\prime_\fm = \alpha_{-1}$. The projection of $\alpha$ onto $\fg_2$
is, effectively, the Maslov form. We will write this projection as $\alpha_{-2}+\alpha_2$ where
$\alpha_{-2}''=0$ and $\alpha_2=\bar\alpha_{-2}$, so that their sum gives the type decomposition of
that projection. We will also fix the metric on $G/G_0$ so that by equation 
\eqref{eq:Maslov} 
\begin{equation}\label{eq:alpha2}
\alpha_{-2}+\alpha_2 = -i\pi\mu\begin{pmatrix} I_2 & 0 \\ 0 & -2\end{pmatrix}.
\end{equation}

From $\alpha$ we can form the \emph{extended Maurer-Cartan form} 
\begin{equation}\label{eq:MCform}
\alpha_\zeta= \zeta^{-2}\alpha_{-2}+\zeta^{-1}\alpha_{-1}+\alpha_0+\zeta\alpha_1+\zeta^2\alpha_2.
\end{equation}
Here $\zeta$ is a rational parameter on the Riemann sphere $\Ci$, so that $\alpha_\zeta$ takes
values in $\fg^\C$ for $\zeta\in\Ct$.
It possesses two symmetries, namely, a real symmetry and $\tau$-equivariance:
\begin{equation}\label{eq:real}
\overline{\alpha_{\zeta}} = \alpha_{\bar\zeta^{-1}},\quad
\tau(\alpha_\zeta) = \alpha_{i\zeta}.
\end{equation}
Consequently we may think of $\alpha_\zeta$ as taking values in
the twisted loop algebra $\Lg$ of $\tau$-equivariant real analytic maps $\xi_\zeta:S^1\to\fg^\C$
possessing the real symmetry. This is a real subalgebra of the complex algebra $\Lgt$ of
$\tau$-equivariant real analytic maps $\xi_\zeta:S^1\to\fg^\C$.
\begin{thm}[\cite{HelR2}]\label{thm:HelR}
Let $\alpha_\zeta$ be the extended Maurer-Cartan form, obtained from any Lagrangian frame,
for a weakly conformal Lagrangian immersion of a surface $f:M\to\CP^2$. Then $f$ is HSL if and only if 
\begin{equation}\label{eq:MC}
d\alpha_\zeta+\frac{1}{2}[\alpha_\zeta\wedge\alpha_\zeta]=0.
\end{equation}
Conversely, suppose $\alpha_\zeta\in\Omega^1_{\tilde M}\otimes \Lg$ has the form \eqref{eq:MCform}
with $\alpha_{-2}''=0$, $\alpha_{-1}\neq 0$ and $\alpha_{-1}''=0$. If $\alpha_\zeta$ satisfies
\eqref{eq:MC} then there exists a \emph{based extended Lagrangian frame} 
\begin{equation}\label{eq:based}
F_\zeta:\tilde M\to C^\omega(S^1,G),\qquad  F_\zeta^{-1}dF_\zeta = \alpha_\zeta,\ 
F_\zeta(0)\in G_0.
\end{equation} 
In that case $F_1$ is the Lagrangian frame for a weakly conformal HSL immersion 
$f:\tilde M\to G/K$.
\end{thm}
The phrase ``weakly conformal'' acknowledges the possibility that $\alpha_{-1}$ (equally $df$, 
and hence the induced metric) may 
vanish at certain points, which we must allow for in the dressing construction: this is explained
in remark \ref{rem:branches} below. Such points will be branch points of the immersion
$f$. 

Notice that since $\alpha_\zeta$ is holomorphic in $\zeta$ throughout $\Ct$, every extended
Lagrangian frame $F_\zeta$ can be extended uniquely holomorphically in $\zeta$ to $\Ct$, in 
which case as a function of $z$ it takes values in $C^\omega(\Ct,G^\C)$. Since $\alpha_\zeta$ is
$\tau$-equivariant and $F_\zeta(0)=I$ we deduce that $F_\zeta$ is also $\tau$-equivariant. 

\section{The spectral data for HSL tori.}

\subsection{Polynomial Killing fields.}
From now on we restrict our attention to the case where $M$ is a 2-torus with conformal class fixed by
$M=\C/\Gamma$, where $\Gamma$ is the period lattice, and fix a conformal coordinate $z$ on $\C$.
In this case HSL immersions of $\C/\Gamma$ into 
$\CP^2$ are of finite type in the sense of integrable systems geometry \cite{HelR3}.
Let us recall what this means. First we need the notion of an adapted polynomial Killing field.
\begin{defn}
An \emph{adapted \pk field} for an extended Maurer-Cartan form
$\alpha_\zeta\in\Omega^1_\C\otimes\Lg$ is a smooth map $\xi_\zeta:\C/\Gamma\to\Lg$ satisfying
\begin{eqnarray}
d\xi_\zeta &  =  & [\xi_\zeta,\alpha_\zeta],\label{eq:Lax}\\
\xi_\zeta  & =  & \zeta^{-4d-2}\alpha_{-2}(\frac{\partial}{\partial z}) + 
\zeta^{-4d-1}\alpha_{-1}(\frac{\partial}{\partial z}) + \ldots,\label{eq:adapted}
\end{eqnarray}
for some non-negative integer $d$. We say $\xi_\zeta$ has \emph{degree} $4d+2$.
\end{defn}
The second condition ensures that the first 
condition is a pair of commuting Lax equations on a finite dimensional subspace of the 
loop algebra, and the arguments of \cite{BurP1} can be adapted to show that this implies 
Hamiltonian integrability (see the appendix \ref{app:dressing} below).
Now we recall that a HSL surface is of \emph{finite type} when it has an extended 
Maurer-Cartan form
which admits an adapted polynomial Killing field. Since this property is gauge invariant for the
gauge group $\caG_0$ it is a property of the Lagrangian angle and not the particular choice of
extended frame.

Adapted \pk fields are not unique: there are infinitely many linearly independent adapted 
\pk fields. 
More generally, if we drop condition \eqref{eq:adapted} we get a large collection
of \pk fields. 
Polynomial Killing fields are the key to encoding the immersion as spectral data. 
The approach we take constructs the spectral data arises from a commutative algebra $\caA$ 
of \pk fields. To understand this algebra we adapt the idea of ``dressing the vacuum 
solution'', as presented in \cite {BurP2}, to take account of quadratic dependence 
of \eqref{eq:MCform} on $\zeta$.

\subsection{Dressing action.}
For any subset $\Delta\subset\Ci$ which is invariant under the real involution
$\zeta\mapsto\bar\zeta^{-1}$ we define
\[
\Lambda^\tau(\Delta,G) = \{g:\Delta\to G^\C:\text{$g$ is analytic},\
\bar{g_\zeta}=g_{\bar\zeta^{-1}},\ \tau(g_\zeta) = g_{i\zeta}\}.
\]
Because the extended Maurer-Cartan form $\alpha_\zeta$ is analytic in $\zeta$ on 
$\Ct$ we can think of every
extended Lagrangian frame $F_\zeta$ as a map from the universal cover $\C$ of our torus into 
$\LHG$. Our dressing action
will be an action of a loop group on the space of HSL immersions of $\C$ into $\CP^2$.
We take the domain to be $\C$ because typically the dressing action does not preserve
periodicity. But dressing does preserve the Maslov form, so that to study dressing orbits of
tori we need only consider the case where $\mu$ is a constant real $1$-form on $\C$ 
(where constant means $\mu'=a\, dz$ for some complex constant $a$).

To each real constant $1$-form $\mu\in\Omega^1_\C$ define $\caE(\mu)$ to be the space of based 
extended Lagrangian frames for HSL immersions of $\C$ into $\CP^2$ with Maslov
form $\mu$, which is equivalent to saying the Maurer-Cartan form of each extended frame satisfies
\eqref{eq:alpha2}. We then define $\caL(\mu)$ to be the space of Lagrangian angles for these 
immersions. One
consequence of theorem \ref{thm:HelR} is that $\caL(\mu)\simeq \caE(\mu)/\caG_0$. The dressing theory
exploits a loop group action on, in our case, a slight enlargement of $\caE(\mu)$ to show that all
HSL tori of finite type and Maslov form $\mu$ lie in the same dressing orbit. The dressing action
is described as follows.

Fix  a real number $0<\e<1$ and define 
\[
C_\e = \{\zeta\in\Ct: |\zeta|=\e\ \textit{or}\ \e^{-1}\}.
\]
This is a union of two circles. The closed annulus in $\Ct$ bounded by $C_\e$ will be called $E_\e$ and 
$I_\e$ will denote the union of closed discs in $\Ci$ with boundary $C_e$. Notice that restriction to
the boundary $C_\e$ allows us to consider both $\LEG$ and $\LIG$ as subgroups of $\LCG$: we shall do
this without further comment. Let $B\subset
G_0^\C$ denote the subgroup of all those matrices in $G_0^\C$ which are upper triangular with positive 
real entries on the diagonal, and define
\[
\LIB = \{g_\zeta\in\LIG:g_0\in B\}.
\]
We recall a fundamental factorisation result concerning the loop 
group $\LCG$. 
\begin{thm}[\cite{McIFW,DorPW}]\label{thm:factor}
The map
\[
\LEG\times\LIB\to\LCG;\quad (g,h)\mapsto gh,
\]
is a diffeomorphism, hence every smooth
$\LCG$-valued function $g(z)$ admits a unique factorisation into the product $g_E(z)g_I(z)$ of smooth
functions with values in $\LEG$ and $\LIB$ respectively.
\end{thm}
As a consequence (see \cite[\S 2]{BurP1})
one obtains the dressing action of $\LIG$, for any $0<\e<1$, on $\LHG$, defined by
\begin{equation}\label{eq:dressing}
\LIG\times\LHG\to \LHG;\quad (g,h)\mapsto g\sharp h = (gh)_E.
\end{equation}
Burstall \& Pedit \cite{BurP1} also show that this action descends to the space $\LHG/G_0$ of cosets
of constant loops.

To apply this to extended Lagrangian frames we first note that, by a standard calculation which
exploits the factorisation theorem (cf.\ equation \eqref{eq:dressalpha} below),
whenever $F_\zeta\in\caE(\mu)$, with Maurer-Cartan form $\alpha_\zeta$ \eqref{eq:MCform}, 
and $g\in\LIG$ the map $g\sharp F_\zeta:\C\to\LHG$ defined by 
\begin{equation}\label{eq:Fdressing}
(g\sharp F_\zeta)(z) = g\sharp (F_\zeta(z))
\end{equation}
has the property that its Maurer-Cartan form has the shape
\begin{equation}\label{eq:degree2}
\zeta^{-2}\alpha_{-2}+\zeta^{-1}\beta_{-1}+\beta_0 +\zeta \beta_1 + \zeta^2 \alpha_2
\end{equation}
for some $1$-forms $\beta_j\in\Omega^1_\C\otimes \fg^\C$ satisfying $\beta_{-1}''=0$ and
$\bar\beta_j = \beta_{-j}$.  In particular, this has the same coefficients of $\zeta^{\pm 2}$ 
as $\alpha_\zeta$ since the adjoint action of $G_0$ on $\fg_2$ is trivial.  Thus $g\sharp F_\zeta$ 
will also be an extended Lagrangian frame for a weakly conformal immersion with the same Maslov form 
provided $\beta_{-1}'$ is not identically zero. This
can fail and in fact $\beta_{-1}= 0$ precisely for what we will call the ``vacuum
solution''. The vacuum solution is the simplest non-trivial
solution of the equations \eqref{eq:MC} for which $\alpha_\zeta$ has the shape \eqref{eq:degree2},
but it does not produce a surface in $G/K$, so there is no corresponding HSL surface. 
For this reason, we must introduce the strict enlargement
$\hat\caE(\mu)$ of $\caE(\mu)$ consisting of those $\LHG$-valued functions $H_\zeta(z)$ on $\C$ which
satisfy \eqref{eq:degree2}, with $\beta_{-1}''=0$ and \eqref{eq:Maslov}. Then \eqref{eq:Fdressing} 
does describe an action of $\LIG$ on
$\hat\caE(\mu)$, and by the remarks above this action descends to an action on the enlarged space
$\hat\caL(\mu) = \hat\caE(\mu)/\caG_0$ 
\begin{defn}\label{defn:vacuum}
The \emph{vacuum solution} in $\hat\caE(\mu)$ is the based extended frame
\begin{equation}\label{eq:Fvacuum}
F^\mu_\zeta = \exp(z\zeta^{-2}A +\bar z\zeta^2 \bar A),\quad
A\, dz = -i\pi\mu'\begin{pmatrix}I_2 & 0 \\ 0 & -2\end{pmatrix}.
\end{equation}
It has Maurer-Cartan form $\alpha^\mu_\zeta = \zeta^{-2} A\, dz +\zeta^2 \bar A\, d\bar z$.
\end{defn}
\begin{rem} 
It is convenient, and makes comparison with \cite{HelR3} easier, for us to write
$-\mu=\bar\mu_0dz+\mu_0d\bar z$, in which case $e^{i(z\bar\mu_0+\bar z\mu_0)}$
is the Lagrangian angle function in the terminology of \cite{HelR3}.
\end{rem}
\begin{lem}\label{lem:vacuum}
Let $F_\zeta\in\hat\caE(\mu)$ be such that its Maurer-Cartan form $\alpha_\zeta$ has $\alpha_{-1}=0$,
then $F_\zeta$ is gauge equivalent to $F^\mu_\zeta$.
\end{lem}
\begin{proof}
By definition of $\hat\caE(\mu)$, $\alpha_{\pm 2}$ satisfy \eqref{eq:Maslov} and therefore
$[\alpha_{-2}\wedge\alpha_2]=0$. The Maurer-Cartan equations \eqref{eq:MC} then imply
\[
d\alpha_0 + \frac{1}{2}[\alpha_0\wedge\alpha_0]=0,
\]
Since we are working over $\C$ this means the term $\alpha_0$ is a pure gauge term, hence
$\alpha_\zeta$ is gauge equivalent to $\alpha^\mu_\zeta$.
\end{proof}
Notice that $F^\mu_1$ takes values in
$\exp(\fg_2)$, which is an $S^1$ subgroup of $K$, hence the vacuum solution corresponds to a constant
map of $\C$ into $\CP^2$. 

We are now in a position to state the dressing theorem for HSL tori in $\CP^2$.
\begin{thm}\label{thm:dressing}
Let $f:\C/\Gamma\to\CP^2$ be a  HSL torus with Maslov form $\mu$. Then it admits a
based extended frame $F_\zeta$ for which, for some $\e >0$ and  $g\in\LIG$, 
$F_\zeta = g\sharp F^\mu_\zeta$.
Consequently, there exists $\chi:\C\to\LIB$ for which
\begin{equation}\label{eq:chi}
gF^\mu_\zeta = F_\zeta\chi.
\end{equation}
\end{thm}
The proof is in appendix \ref{app:dressing}. 
\begin{rem}\label{rem:branches}
The HSL maps $f:\R^2\to\CP^2$ which are constructed by dressing might have branch points. To see
why, recall that branch points occur precisely at the points where $\alpha_{-1}$ vanishes. We can
express this in terms of the coefficients of the
factor $\chi=(I+\chi_1\zeta+\ldots)\chi_0$ in \eqref{eq:chi}. We have
\begin{equation}\label{eq:dressalpha}
\alpha_\zeta = F_\zeta^{-1}dF_\zeta = \Ad\chi_\zeta\cdot\alpha^\mu_\zeta -d\chi.\chi^{-1},
\end{equation}
and so
\[
\alpha_{-1} = \Ad\chi_0\cdot[\chi_1,A].
\]
Now $\ker(\ad A)\cap\fg_{-1}=\{0\}$, so that $\alpha_{-1}$ vanishes precisely where $\chi_1$
vanishes. This is at worst on a proper analytic subset of $\R^2$, but we know of no proof that for
tori it must be non-zero everywhere, so we cannot rule out the possibility of branch points. In the
case of HSL tori in $\R^4$ branch points can also occur \cite{McIR}. However, for minimal
Lagrangian tori in $\CP^2$ branch points cannot occur, since in that case the formula for 
$\alpha_{-1}$ has
the form $\Ad\chi_0\cdot A$ (with the appropriate interpretation of $\chi_0$ and $A$ for that case),
and $\chi_0$ is invertible, hence nowhere zero, on $\R^2$.
\end{rem}
The previous theorem provides the following extremely useful characterisation of \pk fields.
\begin{lem}\label{lem:pkf}
A map $\xi:\C\to\Lgt$ is a \pk field for 
$F_\zeta\in\hat\caE(\mu)$ if and only if $\xi = \Ad\chi\cdot\eta$ for
some $\eta\in\Lgt$ for which $[A,\eta]=0$ and $\Ad\chi\cdot\eta$ is a Laurent polynomial in
$\zeta$ (of bounded degree as $z$ varies).
\end{lem}
The proof is a direct adaptation of lemma A.1 in \cite{McI98}, but the origin of the idea is the
Zakharov-Shabat dressing argument for zero curvature equations. The idea is that
$\eta=\Ad\chi^{-1}\cdot\xi$ satisfies
\[
d\eta + [\eta,\alpha_\zeta^\mu]=0,
\]
and from this it follows that $d\eta=0$ and $[A,\eta]=0$.

\subsection{The algebra generated by \pk fields.}
Our aim is to produce a unital commutative complex matrix algebra $\caA$ which contains all the 
$\tau$-equivariant \pk fields, and ideally $\caA$ should be no larger than necessary. Since the 
Lax equation \eqref{eq:Lax} above is preserved under 
matrix multiplication it is possible to construct $\caA$ so that its elements still satisfy
\eqref{eq:Lax}, but at the expense of losing $\tau$-equivariance. Equation \eqref{eq:Lax} is also
satisfied by any $z$-independent multiple of the identity $I\in\gl_3(\C)$, therefore it is natural
in the first instance to consider \pk fields taking values in the loop algebra
\[
\caR_\e = \{\eta:C_\e\to \gl_3(\C)|\ \eta\ \text{is algebraic}\},
\]
where the terminology ``$\eta$ is algebraic'' means $\eta$ is analytic on $C_\e$ and
extends meromorphically into $I_\e$ where it may only have poles, if any, at $0$ and $\infty$. 
All the $\tau$-equivariant \pk fields take values in $\caR$. Moreover, a simple calculation 
shows that if $\xi_\zeta,\eta_\zeta$ are $\tau$-equivariant \pk fields then 
\[
\xi_{-\zeta}\eta_{-\zeta} = \tau^2(\xi_\zeta\eta_\zeta).
\]
Now we recall that $\tau^2 =\sigma$ \eqref{eq:sigma}, so we may
restrict our attention to \pk fields with values in the subalgebra 
\[
\caR_\e^\sigma = \{\eta_\zeta\in\caR: \sigma(\eta_\zeta) = \eta_{-\zeta}\}.
\]
We will use $\caK^\sigma$ to denote the set of polynomial Killing fields 
with values in $\caR_\e^\sigma$. This is a unital algebra which is independent of $\e$: it contains 
the identity matrix as well as 
all our original \pk fields and their products. The Lax equation \eqref{eq:Lax}, together with the
algebraic properties of $\chi_\zeta$, ensures that $\Ad\chi_\zeta$ maps $\caK^\sigma$ into 
$\caR_\e^\sigma$. As a
consequence of the argument used to prove lemma \ref{lem:pkf}, a map $\xi:\C\to\caR_\e^\sigma$ which
is algebraic in $\zeta$ satisfies
the Lax equation \eqref{eq:Lax} if and only if $\eta=\Ad\chi^{-1}\cdot\xi$ satisfies $d\eta=0$ and
$[\eta,A]=0$. It follows that $\caK^\sigma$ is non-abelian,
because the centraliser $\fz(A)$ of $A$ in $\gl_3(\C)$
is not abelian. To find a maximal abelian subalgebra
$\caA^\sigma\subset\caK^\sigma$ we fix maximal torus $\fh\subset\fz(A)$. It does not matter which
one, but for convenience we choose the maximal torus of diagonal matrices in $\gl_3(\C)$. We let
$\caH_\e^\sigma\subset\caR_\e^\sigma$ denote the elements of $\caR_\e^\sigma$ with values in $\fh$.
Now we can define
\[
\caA^\sigma = \{\xi=\Ad\chi\cdot\eta:\eta\in\caH_\e^\sigma, \xi\ \text{Laurent poly.\ in}\ \zeta\}.
\]
This is a commutative unital $\C$-algebra of \pk fields. 

In fact, $\caA^\sigma$ is best thought of as a complex $1$-parameter family of subalgebras of 
$\caR_\e^\sigma$:
\[
\caA^\sigma(z) = \{\xi(z):\xi\in\caA^\sigma\}\subset\caR_\e^\sigma.
\] 
These are all isomorphic since the Lax equation \eqref{eq:Lax} is equivalent to the condition
\begin{equation}\label{eq:evolution}
\xi_\zeta(z) = \Ad F_\zeta(z)^{-1}\cdot\xi_\zeta(0).
\end{equation}
Moreover, the algebra homomorphism $\caA^\sigma\to\caA^\sigma(0)$ which maps $\xi(z)$ to 
$\xi(0)$ is injective for the same reason, so we may identify each $\caA^\sigma(z)$ with 
$\caA$ itself.

\subsection{Untwisting.}\label{sec:untwist}

The geometric interpretation of $\caA^\sigma$ simplifies by eliminating 
the $\sigma$-equivariance of these \pk fields using ``untwisting'', which is an algebra
isomorphism $\caR_\e^\sigma\simeq\caR_{\e^2}$. This ``partially untwists'' the
$\tau$-equivariant loops, i.e., they retain a symmetry which distinguishes them. Because we need 
to work with the extended frame as well
as the \pk fields we must introduce this untwisting at the loop group level as well as the algebra
level.

First we define
\[
\kappa_\zeta = \begin{pmatrix} \zeta I_2 & 0 \\0 & 1\end{pmatrix},
\]
the significance being that $\kappa_\zeta$ is the $1$-parameter subgroup of $SU(3)$ for which
conjugation by $\kappa_{-1}$ yields the involution $\sigma$. By this property whenever $g_\zeta \in
\LCGsig$ the quantity $\kappa_\zeta^{-1}g_\zeta\kappa_\zeta$ is a function of $\zeta^2$: we will
write this, for the moment, as
\begin{equation}\label{eq:untwist}
\hat g_\lambda = \kappa_\zeta^{-1}g_\zeta\kappa_\zeta,\quad \lambda = \zeta^2.
\end{equation}
At the loop group level, untwisting is the Lie group isomorphism
\[
\LCGsig\to\LCGun;\quad g_\zeta\mapsto \hat g_\lambda.
\]
We are mainly interested in the corresponding Lie algebra isomorphism, which can be extended
to $\caR_\e^\sigma$ and provides the algebra isomorphism
\begin{equation}\label{eq:untwistalg}
\caR_\e^\sigma\to\caR_{\e^2}:\quad \eta_\zeta\mapsto \hat\eta_\lambda =
\kappa_\zeta^{-1}\eta_\zeta\kappa_\zeta.
\end{equation}
The $\tau$-equivariant loops are characterised after untwisting by the symmetry
\begin{eqnarray}
\hat g_{-\lambda}& = & S_\lambda^{-1}\hat g_\lambda^*S_\lambda,\quad g_\zeta\in 
\LCGsig,\label{eq:tau_untwisted}\\
\hat\eta_{-\lambda} & = & -S_\lambda^{-1}\hat\eta_\lambda^t S_\lambda,\quad
\eta_\zeta\in\caR^\sigma,
\end{eqnarray}
where 
\[
S_\lambda = \begin{pmatrix} 0 & i\lambda & 0 \\ -i\lambda & 0 & 0 \\ 0 &0&1\end{pmatrix}.
\]
For each $z\in\C$ we will define $\caA(z)\subset\caR_{\e^2}$ to be the image of $\caA(z)^\sigma$ under 
untwisting: this is independent of $\e$ because its elements are Laurent polynomials. Because of 
equation \eqref{eq:evolution} we also need to consider the untwisted version
of the dressing factorisation \eqref{eq:chi}, which we write as
\begin{equation}\label{eq:dress_fact}
\hat g_\lambda\hat F^\mu_\lambda = \hat F_\lambda\hat \chi_\lambda.
\end{equation}
It is easy to check that
\[
\hat F^\mu_\lambda = \exp(z\lambda^{-1}A +\bar z\lambda\bar A),
\]
and that $g_\zeta$ belongs to $\LIG$ if and only if $\hat g_\lambda$ extends holomorphically into  
the interior of $C_{\e^2}$ and $\hat g_0$ has the shape
\begin{equation}\label{eq:shape}
\begin{pmatrix} a & * &*\\0&a^{-1}&*\\0&0&1\end{pmatrix},\quad a\in\C.
\end{equation}
From now on we will work exclusively in the untwisted setting. For this reason, the ``hat''
notation used above will be dropped.

\subsection{The spectral curve.}

The spectral data is the geometric realisation of the algebra $\caA$ which encodes
the isomorphism class of $\caA$, its symmetries and
the natural gradings $\caA$ carries as an algebra of Laurent polynomials. 
For ease of notation we introduce the algebra homomorphism
\begin{equation}\label{eq:h}
h:\caA\to\caH_{\e^2};\quad \xi\mapsto \Ad\chi^{-1}\cdot\xi.
\end{equation}
We denote the $j$-th diagonal entry of the diagonal matrix $h(\xi)$ by $h_j(\xi)$: this is an
analytic function on $C_{\e^2}$ which extends meromorphically into $I_{\e^2}$. To aid our
discussion we will set $\caC = C^\omega(C_{\e^2},\C)$.

To begin, notice that $\caA$ contains a copy of the algebra $\caB = 
\C[\lambda^{-1},\lambda]$, as scalar multiples of the identity, and the
three prime ideals
\[
\caI_j = \{\xi\in\caA:h_j(\xi)=0\},\quad j=1,2,3, 
\]
which may be trivial. It also possess two involutions $\hat\rho$ and $\hat\tau$, 
the first real linear and the second complex linear, defined by
\begin{equation}\label{eq:rho_tau}
\hat\rho(\xi_\lambda) = (\xi_{\bar\lambda^{-1}})^\dagger,\quad 
\hat\tau(\xi_\lambda) = \Ad S_\lambda^*\cdot(\xi_{-\lambda})^t.
\end{equation}
\begin{lem}\label{lem:rho_tau}
$\hat\rho$ and $\hat\tau$ are algebra automorphisms of $\caA$ (over $\R$ and $\C$
respectively) which commute, and they each preserve the subalgebra $\caB$. While $\hat\rho$ 
preserves $\caI_1,\caI_2,\caI_3$, $\hat\tau$ preserves $\caI_3$ and swaps $\caI_1$ with $\caI_2$.
\end{lem}
\begin{proof}
Both $-\hat\rho$ and $-\hat\tau$ are Lie algebra automorphisms which fix 
the untwisted Maurer-Cartan form $\alpha_\lambda$. Hence if $\xi$ is any (untwisted) \pk field
then so is $-\hat\rho(\xi)$ and $-\hat\tau(\xi)$. Since the Lax equations \eqref{eq:Lax} are
linear in $\xi$ it follows that $\hat\rho(\xi)$ and $\hat\tau(\xi)$ are also \pk fields.
Further, it is easy to show that $S_\lambda^\dagger = S_{\bar\lambda^{-1}}$ and therefore
\[
\hat\rho\hat\tau(\xi_\lambda) = 
\Ad (\bar S_{\bar\lambda^{-1}})^{-1}\cdot \bar\xi_{-\bar\lambda^{-1}}=
\Ad S_\lambda^*\cdot\bar\xi_{-\bar\lambda^{-1}} = \hat\tau\hat\rho(\xi_\lambda).
\]
Next, $\chi_\lambda$ possesses the symmetries
\begin{equation}\label{eq:group_symmetries}
\chi_{\bar\lambda^{-1}} = (\chi_{\lambda}^{-1})^\dagger,\quad
\chi_{-\lambda} = S_\lambda^{-1}\chi_\lambda^* S_\lambda.
\end{equation}
From these it follows that for any $\xi\in\caA$
\[
\hat\rho(\xi) = \Ad\chi\cdot \hat\rho(h(\xi)),\quad
\hat\tau(\xi) = \Ad\chi\cdot \hat\tau(h(\xi)),
\]
where we have used the fact the the definitions of $\hat\rho$ and $\hat\tau$ apply equally well
to $\caH_{\e^2}$. It is easy to see from this that $\hat\rho$ and $\hat\tau$ preserve $\caB$,
since $h(\caB) = \caB$. Clearly $\hat\rho(h(\xi))$ has a zero on its $j$-th diagonal entry if
and only if $h(\xi)$ does (i.e., $h_j(\xi)=0$), therefore $\hat\rho$ preserves $\caI_j$.
Finally, on $\caH_{\e^2}$,
\[
\hat\tau(\begin{pmatrix} a(\lambda) &0&0\\ 0& b(\lambda) &0 \\
0&0&c(\lambda)\end{pmatrix}) = 
\begin{pmatrix} b(-\lambda) &0&0\\0 & a(-\lambda) &0 \\
0&0&c(-\lambda)\end{pmatrix}
\]
and therefore $\hat\tau$ swaps $\caI_1$ with $\caI_2$ but fixes $\caI_3$.
\end{proof}

The geometric realisation is that we have an affine curve $X_A=\Spec(\caA)$ with a natural
rational function $\lambda:X_A\to\C^\times$ (dual to $\caB\subset\caA$) and $X_A$ comes
equipped with a holomorphic involution $\tau$ and a real involution $\rho$, which commute. 
This curve is reducible if any of the ideals $\caI_j$ is non-trivial: the irreducible
components are isomorphic to $\Spec(\caA_j)$. The involution $\tau$ maps $\Spec(\caA_1)$ to
$\Spec(\caA_2)$, which can be either distinct or identical (i.e., either $\caI_1$ and $\caI_2$
are distinct or identical) and $\lambda$ is non-constant on every irreducible
component of $X$.

We will define the spectral curve of $\caA$ to be the completion $X$ of $X_A$ by smooth points. 
Its irreducible component corresponding to the prime ideal $\caI_j$ will be denoted by $X_j$.
Clearly
the points of completion must lie over $\lambda=0,\infty$. The nature
of this completion is determined by the valuations on its field of rational functions $\C(X)$
(equally, the field of fractions of $\caA$) which correspond to measuring the degree of each
rational function at the points of completion. Because $\caA$ need not be an integral domain we 
must take a little care when defining these.

For any algebraic element $a\in \caC$ define its degree about
zero $\deg_0(a)$ to be the degree at $\lambda=0$ of its meromorphic extension into $I_{\e^2}$,
and similarly define its degree about infinity $\deg_\infty(a)$. Using these notions we define
\begin{eqnarray}\label{eq:nu}
\nu^0_j:\caA-\caI_j\to \Z,&\quad& \nu^0_j(\xi) = \deg_0(h_j(\xi)), \\
\nu^\infty_j:\caA-\caI_j\to \Z,&\quad& \nu^\infty_j(\xi) = \deg_\infty(h_j(\xi)).
\end{eqnarray}
It is easy to check that these six functions have the valuation properties:
\begin{enumerate}
\item $\nu(\xi_1\xi_2) = \nu(\xi_1)+\nu(\xi_2)$,
\item $\nu(\xi_1+\xi_2) \geq \min\{\nu(\xi_1),\nu(\xi_2)\}$.
\end{enumerate}
\begin{lem}\label{lem:fractions}
Each of $\nu_j^0$ and $\nu_j^\infty$ induces a valuation on the field of fractions $\caF_j$
of the integral domain $\caA_j$. Further, 
\begin{equation}\label{eq:nu_rho_tau}
\nu_1^0\circ\hat\tau = \nu_2^0,\quad \nu_3^0\circ\hat\tau=\nu_3^0,\quad
\nu_j^0\circ\hat\rho = \nu_j^\infty.
\end{equation}
These valuations correspond to $6$ smooth points on $X$:
$P_j$ corresponding to $\nu_j^0$ and $Q_j$ corresponding to $\nu_j^\infty$. These points have
the properties
\begin{equation}\label{eq:points}
\lambda(P_j)=0,\ \lambda(Q_j)=\infty,\ Q_j=\rho P_j,\ P_2=\tau P_1,\ \tau P_3 = P_3
\end{equation}
Consequently the rational function $\lambda:X\to\Ci$ has degree $3$. 
\end{lem}
\begin{proof}
An element of the field of fractions of $\caA_j$ can be written $[\xi_1]/[\xi_2]$ where
$[\xi_1],[\xi_2]\in\caA_j$ and $[\xi]$ denotes $\xi+\caI_j$. We define 
\[
\nu_j:\caF_j-\{0\}\to\Z;\quad \nu_j([\xi_1]/[\xi_2]) = \nu_j(\xi_1)-\nu_j(\xi_2),
\]
for $\nu_j$ either of $\nu_j^0$ or $\nu_j^\infty$.
It is easy to check that this is well-defined and retains the valuation properties. The
symmetries \eqref{eq:nu_rho_tau} follow from lemma \ref{lem:rho_tau} and the fact that
$\deg_0(a(\lambda)) = \deg_\infty(\bar a(\bar\lambda^{-1}))$ for any algebraic $a\in
C^\omega(C_{\e^2},\C))$. The valuations $\nu_j^0,\nu_j^\infty$ correspond to smooth points
$P_j,Q_j\in X_j$. On $\Ci$ the points $0,\infty$ correspond to the
valuations $\deg_0,\deg_\infty$ on $\C(\lambda)$, which we identify with the field of fractions 
of $\caB$. Now on $\C(\lambda)\subset\caF_j$ we have $\deg_0=\nu_0^j$ and 
$\deg_\infty = \nu_\infty^j$,
so $\lambda(P_j)=0$ while $\lambda(Q_j)=\infty$. The symmetries in \eqref{eq:points} follow from
\eqref{eq:nu_rho_tau}. Finally, matrix multiplication of vectors makes
$\caB\otimes\C^3$ a rank three $\caB$-module. It is also clearly a faithful
$\caA(0)$-module, and therefore $\caA$ is at most rank three as a $\caB$ module. It follows
that $\lambda$ is at most degree three. Since it possesses at least three points over $0$ it must
be exactly degree three. Finally, the symmetries \eqref{eq:points} follow easily from the proof
of lemma \ref{lem:rho_tau}.
\end{proof}
We can also import the argument from \cite[Prop.6]{McI95} to deduce that $X$ can only be
disconnected if $f$ is not linearly full. Since immersed Lagrangian surfaces in $\CP^2$ must be
linearly full, $X$ is connected. In particular, this disposes of the spectral curve for the vacuum
solution, since it is a disjoint union of three copies of the Riemann sphere.
We will use the term ``spectral curve'' to mean the totality of information
$X,\lambda,\tau,\rho$, which includes the points over $\lambda=0$ and $\lambda=\infty$.
Up to this point we have proved the following.
\begin{thm}
The spectral curve $(X,\lambda,\rho,\tau)$ for a HSL torus in $\CP^2$ is a complete
connected algebraic curve $X$ equipped with a degree $3$ rational function $\lambda$, 
a real involution $\rho$
and a holomorphic involution $\tau$. The function $\lambda$ has symmetries
\begin{equation}\label{eq:lambda_symmetries}
\overline{\rho^*\lambda} = \bar\lambda^{-1},\quad \tau^*\lambda = -\lambda.
\end{equation}
Further, $\lambda$ has distinct zeroes (and therefore distinct poles) each of which is a 
smooth point and only one of which is fixed by $\tau$.
\end{thm}

\subsection{The linear family of line bundles.}

While the spectral curve determines $\caA$ it does not carry the information which
distinguishes $\caA(0)$ from $\caA(z)$: this is encoded in the representation of each
$\caA(z)$ on its module $\caM = \caB\otimes\C^3$. The geometric realisation
of this is a rank one torsion free\footnote{In the case that $X$ is reducible ``torsion free''
means the restriction to each irreducible component is torsion-free.} coherent sheaf 
$\caL^z$ over $X$. From now on whenever we use the phrase ``line bundle'' it should be
interpreted to include maximal rank 1 torsion free coherent sheaves whenever $X$ is not smooth. 
The term ``maximal'' means the sheaf is not the direct image of any sheaf on a less singular
curve than $X$, equally, $\Hom(\caL,\caL)$ is trivial. 

To begin, for each $z$ the $\caA(z)$-module $\caM$ determines a line bundle $\caL^z$ over the 
affine curve $X_A$ in the usual way, i.e., the stalk $\caL^z_P$ at $P\in X_A$ is the
localisation $\caM(z)_P$ of $\caM$ at the ideal $\fm(z)_P\in\Spec(\caA)$ corresponding to the
point $P$. Our task is to extend this to the completion $X$ in a natural way. The key is to use
filtrations of $\caM$ which are compatible with the valuations determining the points of
completion. We proceed as follows.

First we introduce the linear map
\begin{equation}\label{eq:s}
s:\caM\to\caC\otimes\C^3;\quad v\mapsto \chi(z)^{-1}v.
\end{equation}
This depends upon $z$, but we will supress this fact in the notation since we will not have a
cause to explicitly refer to this dependence. Now
define $s_j(v)\in\caC$ to mean the $j$-th entry in the vector $s(v)$. This is compatible
with the homomorphism $h$ in \eqref{eq:h}, since $s(\xi v) = h(\xi)s(v)$. Consequently the
subspace
\[
\caS_j(z) = \{v\in\caM:s_j(v)=0\}
\]
is an $\caA$-submodule of $\caM$ with the property that $\caI_j\caM\subset\caS_j(z)$. Therefore
the quotient space $\caM_j(z)=\caM/\caS_j(z)$ is an $\caA_j(z)$-module.

On each of these we define two functions
\begin{eqnarray}\label{eq:mu}
\mu_j^0:\caM_j(z)-\{0\}\to\Z;&\quad &\mu_j^0([v]) = \deg_0(s_j(v)),\\
\mu_j^\infty:\caM_j(z)-\{0\}\to\Z;&\quad &\mu_j^\infty([v]) = \deg_\infty(s_j(v)),
\end{eqnarray}
where $[v]$ stands for $v+\caS_j$. These have the properties
\[
\mu([v] + [w]) \geq \min\{\mu([v]),\mu([w])\},\quad \mu([\xi v]) = \nu([\xi]) + \mu([v]),
\]
whenever $[v],[w]\in\caM_j(z)$ and $[\xi]\in\caA_j(z)$. Because of these properties each
function extends naturally the module of fractions $\tilde\caM_j(z) = 
\caF_j(z)\otimes_{\caA_j(z)}\caM_j(z)$. We define
\begin{eqnarray}\label{eq:LPQ}
\caL^z_{P_j}&  =& \{ v\in\tilde\caM_j(z)-\{0\}:\mu_j^0(v)\geq 0\}\cup\{0\},\\
\caL^z_{Q_j}&  =& \{ v\in\tilde\caM_j(z)-\{0\}:\mu_j^\infty(v)\geq 0\}\cup\{0\}.
\end{eqnarray}
It follows that by attaching these stalks over the appropriate points of $X$ we obtain a
sheaf of $\caO_X$-modules which we define to be $\caL^z$.
\begin{lem}\label{lem:caL}
For each $z_0\in\C$ the sheaf $\caL^{z_0}$ is a maximal, rank one, torsion free coherent sheaf with 
$h^0(\caL^{z_0})=3$ and $h^1(\caL^{z_0})=0$. Thus, by the Riemann-Roch theorem, it has degree 
$\deg(\caL^{z_0})=g+2$, where $g$ is the arithmetic genus of $X$.
\end{lem}
\begin{proof}
From the construction it is clear that $\caL^{z_0}$ is a torsion free coherent sheaf. That it has
rank one follows from the fact that $\caM$ is a rank three $\caB$-module and $\lambda$ has
degree three. To see that it is maximal it suffices to show that the sheaf
$\Hom(\caL^{z_0},\caL^{z_0})$ is trivial. For any affine open $U\subset X$
consisting only of smooth points one knows that $\Hom(\caL^{z_0}_U,\caL^{z_0}_U)\simeq \caO_U$,
(since $\caL^{z_0}_U$ is locally free and rank one) so
since any singularities lie in $X_A$ it suffices to check that this for $U=X_A$, and this
amounts to checking that if $\eta\in\Hom_\caA(z_0)(\caM,\caM)$ then $\eta$ represents
multiplication by an element of $\caA(z_0)$. We use the natural embedding
\[
\Hom_\caA(z_0)(\caM,\caM)\subset \Hom_\caB(\caM,\caM)\simeq \gl_3\otimes\caB
\]
to view $\eta\in\gl_3\otimes\caB$. This must commute with every element of $\caA(z_0)$, but
$\caA(z_0)$ is a maximal abelian subalgebra of $\gl_3\otimes\caB$, hence $\eta\in\caA(z_0)$.

The remainder of the lemma is equivalent to the assertion that the direct image
$\caE^{z_0}=\lambda_*\caL^{z_0}$ is a trivial rank three vector bundle over $\Ci$, since
$h^i(\caE^{z_0})=h^i(\caL^{z_0})$. We will prove this by
showing that if $\sigma_1^{z_0},\sigma_2^{z_0},\sigma_3^{z_0}$ are the three sections of
$\lambda_*\caL^{z_0}$ over $\Ct$ corresponding to the generators $e_1,e_2,e_3$ for
$\caM$ over $\caB$ (the standard basis vectors for $\C^3$) then they extend holomorphically to
$\Ci$ and span $H^0(\caL^{z_0})$. By our construction of $\caL^{z_0}$ over $\lambda=0,\infty$
this means we must consider each $\chi_\lambda^{-1}e_j$ about $\lambda=0,\infty$. But these are
the columns of matrix $\chi_\lambda^{-1}$, so each is holomorphic and
non-vanishing about $\lambda=0,\infty$. Therefore the only linear combinations over $\caB$ of these
which remain holomorphic at both ends are those with constant coefficients. Hence
\begin{equation}\label{eq:Eframe}
H^0(\caL^{z_0})=H^0(\caE^{z_0}) = \Sp_\C\{\sigma_1^{z_0},\sigma_2^{z_0},\sigma_3^{z_0}\}.
\end{equation}
\end{proof}
It will be useful for us to know later how to characterise the basis of sections $\sigma_j^z$. 
\begin{lem}\label{lem:basis}
Up to scaling, the sections $\sigma_1^z$, $\sigma_2^z$ and $\sigma_3^z$ are the 
unique global sections of, respectively, $\caL^z(-P_2-P_3)$, $\caL^z(-P_3-Q_1)$ and 
$\caL^z(-Q_1-Q_2)$.
\end{lem}
\begin{proof}
That these sections vanish at the points indicated follows directly from an inspection of the
columns of $\chi^{-1} e_j$ evaluated at $\lambda=0$ and $\infty$, which can be read off
\eqref{eq:shape} since this is the shape of $\chi_0$ and $\chi_\infty = \bar\chi_0^*$. 
For the same reason, a constant vector $v\in\C^3$ corresponds to a global section vanishing at
$P_2+P_3$ only if $v$ is a scalar multiple of $e_1$. Similarly, the other two vanishing conditions
can only occur when $v=e_2$ and $v=e_3$ respectively.
\end{proof}
As a consequence we identify the direct image $\caE(z)$ with the trivial bundle $\caE$ over $\Ci$
using the frame $\sigma^z$ obtained from \eqref{eq:Eframe}.

Although we constructed $X$ and $\caL^z$ purely algebraically, we can now invoke Serre's GAGA
principle and treat them as complex analytic objects too. This is essential if we wish to
understand the family $\caL^z$, since it is not in general algebraic. This perspective allows
us to interpret the map $s$ as a trivialisation $\varphi^z$ for $\caL^z$ over 
$U=\lambda^{-1}(I_{\e^2})$,
which we may assume without loss of generality contains no branch points of $\lambda$ and is a
disjoint union of open discs, one for each point $P_j$ or $Q_j$. For then,
by definition, $v\in\tilde\caM_j$ represents a holomorphic section over $U$ if and only if
the components $s_j(v)$ of $\chi_\lambda^{-1}v$ are holomorphic
functions on $I_{\e^2}$. Since $s$ is linear over $\caA(z)$ this defines a trivialisation of
$\caL^z$ over $U$. Because $U$ contains no branch points this pushes down to a trivialisation
$\Psi_U^z$ for the direct image $\caE(z)$. 
The transition relation between this and the frame $\sigma^z$ is clearly
\begin{equation}\label{eq:transition1}
\chi(z)\Psi_U^z = \sigma^z,\quad \text{over}\ I_{\e^2}.
\end{equation}
This provides us with a direct link to the dressing factorisation \eqref{eq:dress_fact}, for in the
analytic category we can replace, away from $0,\infty$, the frame $\sigma^z$
by the frame $\Psi_A^z$ corresponding to the columns of $F(z)$, 
i.e., $\Psi_A^z(\sigma_j) = F(z)e_j$, and therefore
\begin{equation}\label{eq:transition2} 
g F^\mu(z)\Psi_U^z = \Psi_A^z,\quad\text{over}\ I_{\e^2}-\{0,\infty\}.
\end{equation}
These transition relations are encoded in the spectral data in the following way. For the
statement of this proposition we assume $\caL$ is a line bundle, for ease of exposition.
\begin{prop}\label{prop:L}
Let $R$ be the ramification divisor of $\lambda$ and write $\caL^z=\caO_X(D)$ for some positive 
divisor $D$. Let $J_\R(X)$ denote the real subgroup $\{L\in\Jac(X):\overline{\rho^*L}\simeq L^{-1}\}$
of $\Jac(X)$. Finally, let $X^s$ be the variety of smooth points on $X$.
\begin{enumerate}
\item We have the linear equivalences
\begin{equation}\label{eq:L_symmetries}
\rho_*D+D\sim R,\quad \tau_*D+D\sim R+P_3-Q_3.
\end{equation}
Further, $\rho$ fixes every point $P\in X$ for which $|\lambda(P)|=1$ and $D$ has no support on
this unit circle.
\item 
The map $\ell:\C\to J_\R(X)$ be given by $\ell(z)= \caL^z\otimes(\caL^0)^{-1}$ is a homomorphism of real
groups. It is completely characterised by
\begin{equation}\label{eq:ltangent}
\frac{\partial \ell}{\partial z}|_{z=0} =
\frac{3i\pi\bar\mu_0}{2}
\left(
\frac{\partial\fA_{P_1}}{\partial\lambda}(0)+
\frac{\partial\fA_{P_2}}{\partial\lambda}(0)
\right)
\end{equation}
Here $\fA_P:X^s\to\Jac(X)$ denotes the Abel map with base point $P$, i.e., 
$\fA_P(Q) = \caO_X(Q-P)$.
\end{enumerate}
\end{prop}
\begin{rem}
(i) To interpret the right hand side of equation \eqref{eq:ltangent} we identify $T_eJ_\R(X)^\C$ with
$T_e\Jac(X)\simeq T_e^{1,0}\Jac(X)$. 
Notice that $\ell_z$ actually takes values in a real subgroup of the Prym variety
$\Prym(X,\tau)\subset\Jac(X)$: since $\tau$ has precisely two fixed points ($P_3$ and $Q_3$) this
Prymian has dimension $g/2$ when $g$ is the arithmetic genus of $X$.\\
(ii) The role of the ramification divisor is to provide a coherent sheaf $\caD=\caO_X(R)$
over $X$ for which
\[
\lambda_*(\Hom_{\caO_X}(\caL,\caD)) = (\lambda_*\caL)^*.
\]
When $X$ is smooth we can take $\caD=\caO_X(R)$ and use the trace map $\C(X)\to\C(\lambda)$ as
part of the dual pairing (see \cite{McI95}). This trace map is essentially the same as the trace in
Serre duality (cf. \cite[II\S 12]{Ser}) so 
when $X$ is singular it may be necessary to interpret $\caO_X(R)$ as $K_X\otimes\lambda^*K_{\Ci}^{-1}$ 
(cf.\ Hartshorne \cite[Ch.IV,\S 2]{Har}), where $K_X$ is the dualising sheaf for $X$.
\end{rem}
\begin{proof}
The equivalence $D+\rho_*D\sim R$ is a direct consequence (see \cite{McI96}) of the real symmetry 
of $\chi$ in \eqref{eq:group_symmetries}, from which we also obtain the conditions on $\rho$. 
In particular, there is a rational function with divisor $D+\rho_*D-R$ which is real
and positive over the unit circle, which is equivalent to the
requirement that $D$ (or a representative in its linear equivalence class) have no points on the
unit circle.

The symmetry $\chi_{-\lambda} = S_\lambda^{-1}\chi_\lambda^*S_\lambda$
gives the second linear equivalence as follows. First, $\chi_\lambda^*$ is the transition function 
for the dual bundle $(\lambda_*\caL)^*$, which is the direct image of $\caO_X(R-D)$ (again, from
\cite{McI96}). Therefore that symmetry means
\[
S_\lambda\tau^*\Psi_U = \Psi_U^*
\]
where $\Psi_U^*$ is the trivialisation of $\caO_X(R-D)$ over $U$ dual to $\Psi_U$. Hence any local
trivialisation of $\caO_X(\tau_*D)$ over $U$ corresponds to a section of $\caO_X(R-D)$ with divisor
$P_1+P_2-Q_1-Q_2$.

For part (ii), we first note that equation \eqref{eq:transition2} tells us that both $\Psi_U^0$ and
$F^\mu(z)\Psi_U^z$ are local trivialisations for $\lambda_*\caL^0$, and since $F^\mu$ is diagonal
both come from local trivialisations of $\caL^0$. Therefore $\caL^z$ is obtained from $\caL^0$ by
twisting it by a line bundle obtained by glueing the trivial bundles over $U$ and $X_A$ together
using the diagonal elements of $F^\mu(z)$ as the transition functions, with the $j$-th diagonal
providing the transition function over the punctured discs about $P_j,Q_j$. This gives a
$1$-cocycle for $\ell(z)$. Differentiating this cocycle gives us
\[
\frac{\partial \ell}{\partial z}|_{z=0} =
\frac{i\pi\bar\mu_0}{2}
\left(
\frac{\partial\fA_{P_1}}{\partial\lambda}(0)
+\frac{\partial\fA_{P_2}}{\partial\lambda}(0)
-2\frac{\partial\fA_{P_3}}{\partial\lambda}(0)
\right).
\]
But $(\partial\fA_{P_j}/\partial\lambda)(0)$ can be identified with the linear map
\[
H^0(\Omega_X)\to\C;\quad \omega\mapsto\Res_{P_j}(\lambda^{-1})
\]
so from the residue theorem, since $P_1+P_2+P_3$ is the divisor of zeroes of $\lambda$, we deduce that
\[
\frac{\partial\fA_{P_1}}{\partial\lambda}(0)
+\frac{\partial\fA_{P_2}}{\partial\lambda}(0)
+\frac{\partial\fA_{P_3}}{\partial\lambda}(0)=0.
\]
Thus we obtain \eqref{eq:ltangent}.
\end{proof}

\subsection{Double periodicity.}\label{sec:double}
As well as the properties of the spectral curve $(X,\lambda)$ stated
above, the double periodicity of $f$ can be encoded into the spectral curve in exactly the same way
as was achieved for minimal Lagrangian tori \cite{McI03}. This involves introducing a particular
singularisation of $X$. Let $O_1,O_2,O_3$ be the three distinct smooth points lying over
$\lambda=1$. Define $X_\fo$ to be the singularisation of $X_\fo$ obtained by identifying together
these points (the singular curve for the modulus $\fo=O_1+O_2+O_3$ \cite{Ser}). For any smooth
point $P\in X_\fo$ let $\caA^\fo_P:X_\fo^s\to\Jac(X_\fo)$ denote the Abel map based at $P$. We can
define a homomorphism of real groups $\ell^\fo:\R^2\to J_\R(X_\fo)$ by insisting that
\begin{equation}\label{eq:lotangent}
\frac{\partial \ell^\fo}{\partial z}|_{z=0} =
\frac{3i\pi\bar\mu_0}{2}
\left(
\frac{\partial\fA^\fo_{P_1}}{\partial\lambda}(0)+
\frac{\partial\fA^\fo_{P_2}}{\partial\lambda}(0)
\right)
\end{equation}
\begin{prop}
Let $(X,\lambda)$ be the spectral curve for a HSL torus $f:\C/\Gamma\to\CP^2$, then $\ell^\fo$ is
$\Gamma$-periodic.
\end{prop}
To prove this we begin by observing that the map $f$ is reconstructed from its spectral data in
exactly the same way as the minimal Lagrangian case (and, more generally, harmonic non-isotropic
tori in $\CP^n$), viz, the line $f(z)\in\CP^2$ is given by $[F_1(z)e_3]$ (modulo an isometry of
$\CP^2$) and therefore, by lemma \ref{lem:basis}  corresponds to the line $H^0(\caL_z(-Q_1-Q_2))$ in
$\P H^0(\caL_z)$ and these projective spaces are identified together by the frame $F_1$.
This leaves only the final identification of $\P H^0(\caL_0)$ with
$\CP^2$ unfixed, so $f$ is only recovered up to an isometry of $\CP^2$.

What we have here is an adaptation of the construction detailed in \cite{McI99} for harmonic
maps (see also \S\ref{sec:explicit} below). It follows that we can apply the tools described 
there to describe our reconstruction of $f$. 
In other words, the HSL torus $f:\C/\Gamma\to\CP^2$ factors through $\ell^\fo$:
\begin{equation}\label{eq:theta}
f:\C/\Gamma\stackrel{\ell^\fo}{\to}\Jac(X_\fo)\stackrel{\theta}{\to}\CP^2,
\end{equation}
where $\theta$ is a certain rational map determined by choosing the divisor $Q_1+Q_2$
to single out the point in $\P H^0(\caL_z)$ determined by the line $H^0(\caL_z(-Q_1-Q_2))$.
More detail is given in \S\ref{sec:explicit} below.

In truth, to establish \eqref{eq:theta} we have to prove \emph{first} that $\ell^\fo$ is 
$\Gamma$-periodic, since \emph{ a priori} it could be that the periodicity is due to $\theta$
rather than $\ell^\fo$. For this we argue as follows. The translated map $f^{w}(z) = f(z+w)$
determines the spectral data $(X,\lambda,\caL^w)$ and therefore $\ell(z)$ is $\Gamma$-periodic. But
$\ell^\fo$ lifts $\ell$ to $\Jac(X_\fo)$ (which is a $(\C^\times)^2$-bundle over $\Jac(X)$)
and from \cite[Lemma 1]{McI99} the holonomy of this lift gives a homorphism $h:\Gamma\to PU(3)$ for
which $f(z+\gamma) = h(\gamma)f(z)$ for all $z\in\C$ and each $\gamma\in\Gamma$. Since
$f(z+\gamma)=f(z)$ and $f$ is linearly full, $h(\gamma)$ is trivial for each $\gamma$, 
whence $\ell^\fo$ is $\Gamma$-periodic.

\section{The reconstruction from spectral data.}

In this section we will complete the proof of theorem \ref{thm:main} by showing that we can start
with spectral data and construct a HSL torus $f$ in a manner
which reverses the procedure given above for extracting the spectral data from
$f$. We finish with some brief comments about the consequences of this classification for the
study of the moduli space of HSL tori. 

\subsection{Reconstruction.}\label{sec:Rec}
Suppose we are given a triple of HSL spectral data $(X,\lambda,\caL)$ as per definition
\ref{defn:spectral_data}. Let us first show that provided $(X,\lambda)$ satisfy conditions 
(a) and (b)
of that definition we can always find a line bundle $\caL=\caO_X(D)$ satisfying condition (d). 
First let us rewrite equations \eqref{eq:lequiv} in terms of the canonical class $K_X$ of $X$.
Since $R\sim K_X+\sum_{j=1}^3(P_j+Q_j)$ a degree $g+2$ divisor $D>0$ satisfies \eqref{eq:lequiv} if
and only if $D'=D-Q_1-Q_2$ satisfies
\begin{equation}\label{eq:leqcan}
D'+\rho_*D'\sim K_X+P_3+Q_3,\qquad (\rho\tau)_*D'\sim D'.
\end{equation}
Further, since \eqref{eq:lequiv} implies $\dim H^0(\caL(-Q_1-Q_2)=1$ we may assume $D'>0$. Notice
that $\deg(D')=g$ so that these equations reduce to those of Sharipov \cite{Sha} for minimal
Lagrangian tori.
\begin{rem}
We take this opportunity to point out an error in \cite[Lemma 2]{McI03}. The expression
for equation (15) in that reference should be 
\[
\mu_*\caL_z\simeq\caL_z^{-1}\otimes\caO_X(R+P_\infty-P_0), 
\]
which, given the conventions of that paper, agrees with \cite{Sha}. 
\end{rem}
\begin{lem}
Whenever $(X,\lambda)$ satisfy conditions (a) and (b) of definition \ref{defn:spectral_data} there
exists a degree $g$ divisor $D'>0$ satisfying \eqref{eq:leqcan}.
\end{lem}
\begin{proof}
Set $\iota=\rho\tau$ and let $V = \{k\in\C(X):(k)>-K_X-P_3-Q_3\}$. Since the class $K_X+P_3+Q_3$ 
is both $\rho$ and $\tau$ invariant $V$ admits two real linear involutions, 
$\bar\rho(k)=\overline{\rho^*k}$ and $\bar\iota(k)=\overline{\iota^*k}$. Since $g$ is even 
there exists a
non-zero $k\in V$ with divisor of zeroes of the form $D'+\rho_*D'$ where $\iota_*D'=D'$ if and only 
if $k$ has a divisor of zeroes which is both $\rho$ and $\iota$ invariant. 
The fixed points of $\bar\rho$ form a real subspace of $V$
of dimension $g$ which is preserved by the involution $\bar\iota$, so there exists a non-zero 
$k\in V$ for which $\overline{\rho^*k}=k$ and $\overline{\iota^*k}=\pm k$. 
\end{proof}
In light of \S\ref{sec:double} above we can give the details
of the double periodicity condition (c) in definition \ref{defn:spectral_data}. With the 
notation of \S\ref{sec:double} let 
$\ell^\fo:\R^2\to\Jac_\R(X_\fo)$ be the unique 
homomorphism of real groups satisfying
\begin{equation}\label{eq:lotangent2}
\frac{\partial\ell^\fo}{\partial w}|_{w=0}=\frac{3\pi i}{2}\left(
\frac{\partial\fA^\fo_{P_1}}{\partial\lambda}(0)+
\frac{\partial\fA^\fo_{P_2}}{\partial\lambda}(0)
\right).
\end{equation}
This equation is to be read in the same sense as \eqref{eq:ltangent}. 
The periodicity condition (c) in definition \ref{defn:spectral_data} is that
$\ell^\fo$ must be doubly periodic.

A comparison between the equations \eqref{eq:lotangent2} and \eqref{eq:lotangent} shows us how to
obtain the Maslov form (and hence the mean curvature).
If $\Gamma$ denotes the period lattice for $\ell^\fo$ then \eqref{eq:lotangent2} equips 
$\R^2/\Gamma$ with
the harmonic $1$-form $dw+d\bar w$ (pushed down to $\R^2/\Gamma$, where it is no longer exact).
This will give the negative of the Maslov form $\mu$.

Now we may state and prove the main result of this section, and so complete the proof of theorem
\ref{thm:main}.
\begin{prop}\label{prop:spectral}
To each triple $(X,\lambda,\caL)$ of HSL spectral data we can canonically assign a HSL torus
$f:\C/\Gamma\to \CP^2$, uniquely up to base point preserving isometries of $\CP^2$, 
in a manner which reverses the construction of the spectral data from $f$.
\end{prop}
Although \eqref{eq:theta} provides a direct construction of the HSL torus from its spectral
data, we could not see a direct way of verifying that this map is HSL. So instead
we will prove this proposition by constructing a coset $gZ^\e_A$ in 
$\LIB/Z^\e_A$ which corresponds 
to a point in the dressing orbit $\caO^\e(\mu)$ and yields the Lagrangian angle $\Omega$ of $f$. 
Because of the way we have specified equation \eqref{eq:lotangent2} we work with a normalised
extended frame for the vacuum solution,
\[
F^w_\zeta = \exp[i\pi(w\zeta^{-2}+\bar w\bar\zeta^2)\begin{pmatrix} I_2 & 0 \\ 0 & -2\end{pmatrix}].
\]
Once we have constructed a coset $gZ^\e_A$ we obtain an extended frame $F_\zeta(w) =
g\sharp F^w_\zeta$ for a Lagrangian angle $\Omega:\C/\Gamma\to G/G_0$. The Maslov form of the 
corresponding HSL torus is then $\mu=-(dw+d\bar w)$ as a harmonic $1$-form on $\C/\Gamma$.
\begin{proof}
We begin by choosing $\e>0$ so that the union of open discs
\[
I_{\e^2} = \{\lambda\in\Ci:|\lambda|<\e^2\text{ or } |\lambda|>\e^{-2}\}
\]
contains no branch points of $\lambda$. Let $\caE$ denote $\lambda_*\caL$. The reality conditions
on $\caL$ are such that $\caE$ is trivial \cite{McI96}: our aim is to produce a global trivialisation 
$\Psi$ and another trivialisation $\Psi_I$ over $I_{\e^2}$ so that the transition relation
\begin{equation}\label{eq:transition}
g_\lambda\Psi_I=\Psi,\quad \lambda\in I_{\e^2}-\{0,\infty\},
\end{equation}
furnishes the coset $gZ^\e_A$ in untwisted form.

From \cite{McI96} we know that, as a consequence of the reality conditions, $H^0(\caL)=H^0(\caE)$ 
possesses a Hermitian inner product $h(\sigma_1,\sigma_2)$
derived from the trace map $\Tr:\caO_X(R)\to\caO_{\Ci}$ \cite{McI95,McI99}. This inner
product is determined only up to a positive real scaling factor, but this will prove to be
sufficient\footnote{This lack of uniqueness reflects the fact that throughout the whole discussion
of the HSL condition we have only used the fact that the metric on $\CP^2$ is $PU(3)$-invariant, 
which only fixes the metric up to scale.}.
Further,
the subspaces $H^0(\caL(-P_2-P_3))$ and $H^0(\caL(-Q_1-Q_2))$ are one dimensional (since the
triviality of $\caE$ means neither $\caL(-\sum P_j)$ nor $\caL(-\sum Q_j)$ have non-trivial global
sections). Moreover, one knows that for $\sigma_1,\sigma_2\in H^0(\caL)$ their inner product
$h(\sigma_1,\sigma_2)$ is zero if $\sigma_1\overline{\rho_*\sigma_2}$ vanishes on any divisor
given by a fibre of $\lambda:X\to\Ci$. Using this property it is straightforward to check that 
$H^0(\caL(-P_2-P_3))$, $H^0(\caL(-P_3-Q_1))$ and $H^0(\caL(-Q_1-Q_2))$ are mutually perpendicular
and we can choose an $h$-unitary frame $\sigma_1$, $\sigma_2$ and $\sigma_3$ from those lines in
$H^0(\caL)$. These are fixed only up to common $\R^+$ scaling and individual $S^1$ rotation. 
Together they provide a global trivialisation 
\[
\Psi:\caE\to \caO_{\Ci}^3;\quad \Psi(\sigma_j)=e_j.
\]
Moreover, because $h$ is non-degenerate we have a canonical isomorphism
$\overline{\rho^*\caE}\simeq\caE^*$, where we are using $\rho$ to represent $\lambda\mapsto
\bar\lambda^{-1}$ on $\Ci$. Under this isomorphism $\overline{\rho^*\Psi}$ is
identified with $\Psi^*$, the trivialisation dual to $\Psi$. 

Now set $U=\lambda^{-1}(I_{\e^2})$. This is a disjoint union of three copies of $I_{\e^2}$, since
$\lambda$ is locally biholomorphic about each pair $P_j,Q_j$. Therefore
any trivialisation $\psi$ for $\caL$ over $U$  induces a 
trivialisation $\Psi_I$ for $\caE$ using the identification 
\[
\caE_{I_{\e^2}}=\caL_U\simeq(\caO_{I_{\e^2}})^3.
\]
Further, since we can change $\psi$ independently over each of the connected components of $U$, it
can be chosen so that $\overline{\rho^*\Psi_I}=\Psi_I^*$. Now consider the transition matrix
$g_\lambda$ in \eqref{eq:transition}.
It satisfies $\overline{g_{\bar\lambda^{-1}}}=g_\lambda^*$ because of the
reality conditions on $\Psi$ and $\Psi_I$. Its value $g_0$ at $\lambda=0$ is upper triangular
because of the zeroes of $\sigma_1$ and $\sigma_2$. Moreover we can adjust $\psi$ over $\lambda=0$
so that $g_0$ has the shape \eqref{eq:shape} with positive diagonal entries. 
Our final task is to show that, by possibly placing further conditions on $\Psi$ and $\Psi_I$, we
also ensure the appropriate $\tau$-symmetry \eqref{eq:tau_untwisted} of $g_\lambda$.

To this end set $E=P_1+P_2-Q_1-Q_2$ and define
\[
\theta_1 = i\lambda^{-1}\sigma_2,\ \theta_2=-i\lambda^{-1}\sigma_1,\ \theta_3=\sigma_3,
\]
which we may think of as globally holomorphic sections of $\caF=\lambda_*\caL(E)$. Indeed, these
globally frame $\caF$. Clearly
\[
\Psi(\theta_1) = i\lambda^{-1}e_2,\ \Psi(\theta_2) = -i\lambda^{-1}e_1,\ \Psi(\theta_3)=e_3,
\]
and therefore $(\tau^*S_\lambda)\Psi(\theta_j)=e_j$, i.e., $(\tau^*S_\lambda)\Psi$ globally
trivialises $\caF$. By assumption
\[
\tau_*\caL(E)\simeq \caO_X(R)\otimes\caL^{-1}
\]
and therefore $\tau^*\caF\simeq\caE^*$. We fix this isomorphism as
\begin{equation}\label{eq:iso}
\tau^*\caF\to\caE^*;\quad S_\lambda\tau^*\Psi\mapsto \Psi^*.
\end{equation}
Clearly $S_\lambda\tau^*\Psi_I$ is a trivialisation for $\tau^*\caF$ over $I$. After our choices
above, we still have freedom to scale $\psi$ independently over $P_1$ and $P_2$ in such a way that
$S_\lambda\tau^*\Psi_I$ maps to $\Psi^*$ under \eqref{eq:iso}. Therefore from $g_\lambda\Psi_I=\Psi$ we
deduce
\begin{eqnarray*}
\tau^*g_\lambda\tau^*\Psi_I&=&\tau^*\Psi\\
\implies (S_\lambda\tau^*g_\lambda S_\lambda^{-1})S_\lambda\tau^*\Psi_I& = &S_\lambda\tau^*\Psi,\\
\implies (S_\lambda\tau^*g_\lambda S_\lambda^{-1})\Psi_I^* & = & \Psi^*.
\end{eqnarray*}
Hence $S_\lambda\tau^*g_\lambda S_\lambda^{-1}=g_\lambda^*$.
\end{proof}

\subsection{Explicit reconstruction.}\label{sec:explicit}
The method by which we have proved proposition \ref{prop:spectral} amounts to an algorithm for 
performing the dressing construction of theorem \ref{thm:dressing} but, as mentioned above, to 
construct the HSL torus this can be circumvented using results from \cite{McI99}, which tell us
that $f:\C/\Gamma\to\CP^2$ must factor through the generalised Jacobian $\Jac(X_\fo)$ in the form
\eqref{eq:theta}. Recall that $\fo=O_1+O_2+O_3$ is the divisor of points lying over $\lambda=1$ and
$X_\fo$ is the singular obtained by identifying these three points together.
We will summarise here what this means in the case where $X$ is smooth: the
description for the case where $X$ may be singular simply requires more technical details.

First we describe the rational map $\theta:\Jac(X_\fo)\to\CP^2$.
The general principle, regardless of whether or not $X$ is smooth, is that $\theta$ assigns to each
$\tilde L$ in an open subvariety $\caU\subset \Jac(X_\fo)$ the point in $\CP^2$ corresponding to 
the line 
\[
H^0(\caL(-Q_1-Q_2)\otimes L)\subset H^0(\caL\otimes L),
\] 
where $L=\pi(\tilde L)\in\Jac(X)$ for
the fibration $\pi:\Jac(X_\fo)\to\Jac(X)$. The open set $\caU$ is the set on 
which $H^0(\caL(-Q_1-Q_2)\otimes L)$ has dimension $1$ and hence $\P H^0(\caL\otimes L)\simeq \CP^2$:
it precisely the extra data carried by $\tilde L$ which determines this isomorphism once it has
been fixed for $\P H^0(\caL)$. When $X$ is smooth $\caU$ is the preimage 
of the complement of a translate of the $\Theta$-divisor on $\Jac(X)$. It was shown in 
\cite{McI99} that this map is 
algebraic, hence rational, and that the real subgroup 
\[
J_\R(X_\fo)=\{\tilde L\in \Jac(X_\fo):\overline{\rho^*\tilde L}\simeq \tilde L^{-1}\}
\]
lies inside $\caU$. Thus $\theta$ is real analytic on $J_\R(X_\fo)$, so that 
$f=\theta\circ\ell^\fo$ is also real analytic. The reality conditions on $\caL$ and $\tilde L$ 
are necessary to identify each $\P H^0(\caL\otimes L)$ with $\CP^2$ as a Hermitian symmetric space.

Given a smooth HSL spectral curve $(X,\lambda)$ we can write $\theta$ in terms of the Riemann
$\Theta$-function for $X$. We first identify
\begin{equation}
\label{eq:Jac}
\Jac(X_\fo)\cong H^0(\Omega_X(\fo))^*/H_1(X-\fo,\Z)\simeq
\C^{g+n}/\Lambda^\prime,
\end{equation}
where $\Omega_X(\fo)$ is the sheaf of mermorphic differentials on $X$ with divisor of
poles no worse than $\fo$ and $\Lambda^\prime$ is a lattice on
$2g+2$ generators. We choose coordinates so that
$\pi:\Jac(X_\fo)\to\Jac(X)$ is covered by the map
\[
\pi:\C^{g+n}\to\C^g;\quad \tilde W=(w_1,\ldots,w_{g+2})\mapsto W=(w_1,\ldots,w_g).
\]
These coordinates amount to making a choice of basis for $H^0(\Omega_X(\fo))$ consisting of a
basis $\omega_1,\ldots,\omega_g$ for $H^0(\Omega_X)$ and two meromorphic differentials
of the third kind $\omega_{g+j}\in H^0(\Omega_X(O_1+O_{j+1}))$. These can be normalised in the
standard manner with respect to a choice of ``a-cycles'' for the punctured curve
$X-\fo$, i.e., so that 
\[
\int_{a_i}\omega_j=\delta_{ij},\ i,j=1,\ldots,g+2.
\]
The real symmetry of $X$ is such that the cycles and differentials can be chosen to allow
$\overline{\rho^*\omega_j} =-\omega_j$ and 
$J_R(X_\fo)\simeq \R^{g+2}/(\Lambda'\cap\R^{g+2})$. It is shown in \cite[p432]{McI99} that this
is a real torus.

Let $\Theta$ be the classical Riemann $\theta$-function on $\C^g$ corresponding 
to the induced isomorphism $\Jac(X)\simeq \C^g/\pi(\Lambda^\prime)$. 
Now let us define $\phi_0(\tilde W) = \Theta(W)$ and 
for $j=1,2$ define
\[
\phi_j(\tilde W) = \exp(2\pi iw_{g+j})\Theta(W+\fA_{O_1}(O_{j+1})).
\]
Each $\phi_j$ represents a global holomorphic section of a line bundle over $\Jac(X_\fo)$, namely, 
the pullback by $\pi$ of (an appropriate translate of) the $\Theta$-line bundle over $J(X)$ 
\cite{McI99}. 

Now let $D'$ be the unique positive divisor in the  class $\caL(-Q_1-Q_2)$. Let
$\kappa\in\C^g$ be the appropriate translation for which $\Theta(\fA_{O_1}(P)+\kappa)$ has divisor of
zeroes $D'$. We set $\tilde\kappa=(\kappa_1,\ldots,\kappa_g,0,0)\in\C^{g+2}$.

As a consequence of the reality condition on $\caL$ in \eqref{eq:lequiv}, the set
\[
\caD(\caL)=\{\text{divisors $D$ on}\ X-\fo:\caO_X(D)\simeq\caL,\ D>0,\ D+\rho_*D\sim_\fo R\}
\]
is non-empty, where ``$\sim_\fo$'' means linear equivalence on $X_\fo$, i.e., 
$D+\rho_*D-R$ is the divisor of a rational function taking the value $1$ at each
point $O_1,O_2,O_3$. In fact $\caD(\caL)$ is identifiable with
\[
\{\tilde\caL\in\Pic_{g+2}(X_\fo):\pi(\tilde\caL) = \caL,\ 
\tilde\caL\otimes\overline{\rho^*\tilde\caL}\simeq\caO_{X_\fo}(R)\},
\]
which is a real slice of the fibre of $\pi:\Pic(X_\fo)\to\Pic(X)$ over $\caL$.
The kernel of the group homomorphism $\pi:J_\R(X_\fo)\to J_\R(X)$ acts freely and 
transitively on this set, hence  $\caD(\caL)\simeq S^1\times S^1$, 

Choose some $D\in\caD(\caL)$. Since $D\sim D'+Q_1+Q_2$ 
there is a rational function $k$ on $X$ with divisor $D-(D'+Q_1+Q_2)$
and $k(P)\Theta(\fA_{O_1}(P)+\kappa)$ has divisor
$D-Q_1-Q_2$. Notice it does not vanish at any $O_j$. Now we define
\begin{eqnarray*}
\theta:\Jac(X_\fo)&\to&\CP^2\\
\theta(\tilde W\bmod\Lambda')& = &
[c_0\phi_0(\tilde W+\tilde\kappa),c_1\phi_1(\tilde W+\tilde\kappa),c_2\phi_2(\tilde W+\tilde\kappa)],
\end{eqnarray*}
where the constants $c_j$ are given by
\begin{equation}\label{eq:c_j}
c_j = \frac{1}{k(O_{j+1})\Theta(\fA_{O_1}(O_{j+1}) +\kappa)},\qquad j=0,1,2.
\end{equation}
Finally, the HSL torus $f:\C/\Gamma\to\CP^2$ corresponding to $(X,\lambda,\caL)$ is obtained 
as the composition $f = \theta\circ\ell^\fo$. Different choices of $D\in\caD(\caL)$ alter the
constants $c_j$ by unimodular multipliers, so that $f$ is determined up to an 
isometry of $\CP^2$ by the data $(X,\lambda,\caL)$.

The explicit formula for the real homomorphism
$\ell^\fo:\C/\Gamma\to\Jac(X_\fo)$ can be easily calculated from \eqref{eq:lotangent2} using the
standard observation that
\[
\frac{\partial\caA_{P_i}}{\partial\lambda}(0) = 
(\Res_{P_i}\lambda^{-1}\omega_1,\ldots,\Res_{P_i}\lambda^{-1}\omega_{g+2})\in\C^{g+2}.
\]
\begin{exam}\label{exam:g=0}  
To illustrate this we will compute the HSL tori in $\CP^2$ which arise from the choice
$X=\Ci$. It is a consequence of \cite[\S 4.2]{McI99} that these will be $S^1\times S^1$-equivariant
(called ``homogeneous'' in \cite{HelR3}). Recall that this means $f:\C/\Gamma\to\CP^2$ 
possesses a real homomorphism $h:\C/\Gamma\to SU(3)$ for which $f(p)=h(p)f(0)$ for all $p\in\C/\Gamma$. 
All such HSL tori were constructed explicitly (in non-conformal coordinates) in \cite[\S 5]{HelR3}. 

Let $\zeta$ denote the natural rational parameter on $\Ci$ and, for a fixed
$a\in\C$ with $0<|a|<1$, set
\[
\lambda = \zeta\frac{(\zeta^2-a^2)}{(\bar a^2\zeta^2-1)}\frac{(\bar a^2-1)}{(1-a^2)}.
\]
This choice ensures that the involutions $\rho^*\zeta = \bar\zeta^{-1}$ and $\tau^*\zeta = -\zeta$ 
have the correct effects on $\lambda$ and that we can take $O_1$ to be $\zeta=1$. The other two
points $O_2,O_3$ over $\lambda=1$ are the two roots of the quadratic
\[
\zeta^2+\frac{1-|a|^4}{1-\bar a^2}\zeta + \frac{1-a^2}{1-\bar a^2}.
\]
We also take $\zeta(P_1)=a$, $\zeta(P_2)=-a$ and $\zeta(P_3)=0$ and so forth for the points $Q_j$. 

The first homology of $\Ci-\fo$ is generated by positively oriented cycles $a_1,a_2$ encircling
$O_2$ and $O_3$ respectively and the dual basis of differentials is given by
\[
\omega_j = \frac{1}{2\pi i}(\frac{1}{\zeta-O_{j+1}}-\frac{1}{\zeta-O_1})d\zeta,\quad
j=1,2.
\]
Thus $\Jac(X_\fo)\simeq\C^2/\Z^2$ and $J_\R(X_\fo)\simeq \R^2/\Z^2$. We can take $\Theta=1$ and
thus
\[
\phi_0(W_1,W_2) = 1,\quad \phi_1(W_1,W_2)=\exp(2\pi iW_1),\quad \phi_2(W_1,W_2)=\exp(2\pi iW_2).
\]
The ramification divisor $R$ is the four point divisor of zeroes of
$d\lambda/d\zeta$, and is given by the roots of
\[
(\zeta^2-C_+)(\zeta^2-C_-),\qquad C_{\pm} = \frac{1}{2\bar a^2}(3-|a|^4\pm\sqrt{(|a|^4-1)(|a|^4-9)}.
\]
Let $R_1+R_2$ be the roots of $\zeta^2-C_+$, then $R=R_1+R_2+\rho_*R_1+\rho_*R_2$. 
Since the genus is zero we can choose $\caL$ to be any degree $2$ line bundle and choose
$D=R_1+R_2$. The rational function $k$ with divisor $R_1+R_2-Q_1-Q_2$ and normalised by $k(O_1)=1$
is given by
\[
k(\zeta) = \frac{(\zeta^2-C_+)}{(\bar a^2\zeta^2-1)}\frac{(\bar a^2 -1)}{(1-C_+)}
\]
Finally, define
\begin{eqnarray*}
U_1=\frac{\partial\caA_{P_1}}{\partial\lambda}(0)& = 
&\frac{1}{2\pi i}\frac{(1+a)(1-|a|^4)}{2a^2(1-\bar a^2)}
\left(\frac{O_2-1}{a-O_2},\frac{O_3-1}{a-O_3}\right),\\
U_2=\frac{\partial\caA_{P_2}}{\partial\lambda}(0)& = 
&\frac{1}{2\pi i}\frac{(1-a)(1-|a|^4)}{2a^2(1-\bar a^2)}
\left(\frac{O_2-1}{a+O_2},\frac{O_3-1}{a+O_3}\right),
\end{eqnarray*}
and set $U = (3\pi i/2)(U_1+U_2)\in\C^2$.
Then there is a maximal lattice $\Gamma\subset\C$ for which
\[
\ell^\fo:\C/\Gamma\to \R^2/\Z^2\simeq J_\R(X_\fo);\quad w\bmod\Gamma\mapsto wU+\bar w\bar U\bmod\Z^2,
\]
and the HSL torus is given by
\begin{eqnarray}
f:\C/\Gamma&\to&\CP^2;\label{eq:homog}\\
f(w)& = &[1,c_1\exp(wA_1-\bar w \bar A_1),c_2\exp(wA_2- \bar w\bar A_2)] \notag
\end{eqnarray}
where $c_j=1/k(O_{j+1})$ and 
\[
A_j = \frac{3\pi i(O_{j+1}^2-1)(1-|a|^4)}{4a(a^2-O_{j+1}^2)(1-\bar a^2)}.
\]
This gives a two real parameter family of conformally embedded homogeneous HSL tori, with 
parameter $a$. To verify that $f$ is indeed HSL, let us point out that \eqref{eq:homog} is certainly a
homogeneous immersion of a torus, and the condition that a map of this form is both conformal and
Lagrangian can be shown to be 
\begin{equation}
A_1^2|c_1|^2 + A_2^2|c_2|^2 + |c_1c_2|^2(A_2-A_1)^2=0.
\end{equation}
Numerical calculations verify that this holds for the quantities above. The HSL condition
follows since every homogeneous conformal Lagrangian torus has harmonic Maslov form.

Notice that the computation above must exclude $a=0$, the minimal
Lagrangian limit, because of the expression for $C_{\pm}$, $U_1$ and $U_2$.
However, when $a=0$ we can set $R_1+R_2= 2.\infty$ (since $Q_1=Q_2=Q_3=\infty$ 
is the ramification point over $\lambda=\infty$ in this case) and take
$U=(\partial\fA_{P_1}/\partial\zeta)(0)$. In this case the calculation simplifies greatly,
and  does indeed produce the unique (up to isometries) homogeneous minimal Lagrangian torus in $\CP^2$.
\end{exam}

\subsection{Brief comments on moduli.}\label{sec:moduli}
The HSL spectral data given in definition \ref{defn:spectral_data} looks very similar
to that for a minimal Lagrangian torus. Indeed, the conditions (a), (b) and (d) reproduce the
conditions for a minimal Lagrangian torus \cite{McI03,Sha} when we force $P_1=P_2=P_3$ (and hence
$Q_1=Q_2=Q_3$). In that limit $\lambda$ is branched over $0$ and $\infty$, so it is no longer a
local parameter at those points. Therefore in this limit \eqref{eq:lotangent} must be interpreted
as a statement about tangent planes in the generalised Jacobian. 
Nevertheless, one can think of the passage from minimal
Lagrangian to HSL as a trifurcation of the zeroes of $\lambda$. This 
allow us to count the expected dimension of the moduli space of spectral curve pairs
$(X,\lambda)$ which admit a HSL torus. For minimal Lagrangian tori the expected dimension is zero 
(so that generically no continuous deformations exist, see
\cite{McI03,CarM}). For HSL tori there are two free real parameters in the count, so we can expect 
each HSL spectral curve pair $(X,\lambda)$ to be able to be deformed in two parameter families. 
This is certainly consistent with what we know for the homogeneous tori in example \ref{exam:g=0}
above.

The spectral genus $g$ also gives us a measure of the deformation space for each HSL torus, for we
are able to move $\caL$ continuously without breaking the double periodicity (indeed, this even
fixes the period lattice $\Gamma\subset\C$). Just as with the minimal Lagrangian case,
the symmetry restrictions on $\caL$ oblige it to lie on a real slice of a translate of the Prym
variety $\Prym(X,\tau)$. Since $\tau$ has exactly two fixed points $g$ is even and this Prymian has
complex dimension $g/2$. This contributes $g/2-2$ real dimensions to the deformation space of a HSL
torus of spectral genus $g$, since we must remove the $2$-parameters corresponding to 
translations of the base point. Altogether this predicts smooth $g/2$-dimensional families of
HSL tori, fibered by $g/2-2$-dimensional leaves (each of which will be a torus) consisting of
HSL tori of the same conformal type. 

\appendix

\section{The dressing theory.}\label{app:dressing}

To adapt the dressing theory from \cite{BurP1} to prove theorem \ref{thm:dressing}
we must do two things. First, adapt
the ``Symes formula'' and second describe the dressing orbit of the vacuum solution
\eqref{eq:Fvacuum}. Throughout this appendix we work with the Lie algebra $\Lcg$ of $\LCG$ and its
subalgebras $\Leg$, $\Lig$ and $\Lib$ (the Lie algebras of, respectively, $\LEG$, $\LIG$ and $\LIB$),
where $\fb$ is the Lie algebra for $B\subset G^\C_0$ defined earlier. 
The Lie algebra version of the group factorisation \eqref{thm:factor} is the Lie algebra direct sum 
decomposition
\begin{equation}\label{eq:Lcg}
\Lcg = \Leg \oplus \Lib;\quad \xi = \xi_E+\xi_I.
\end{equation}
Although elements on $\Lcg$ are pairs of loops, the reality condition allows us to identify each one
uniquely with a loop on the radius $\e$ circle about $\zeta=0$. We will do this frequently, and in
particular for any $\xi\in\Lcg$ the notation $\zeta^k\xi$ will mean the element of $\Lcg$ obtained
by multiplying the loop on the radius $\e$ circle by $\zeta^k$. Therefore the set
\[
\Ltg = \{\xi\in\Lcg:\zeta^2\xi\in\Lig\} 
\]
describes a vector subspace of $\Lcg$.
We let $\hat\caE$ stand for $\cup_{\mu}\hat\caE(\mu)$ and let $\hat\caL$ stand for $\hat\caE/\caG_0$.
We denote by $[F_\zeta]\in\hat\caL$ the gauge equivalence class $F_\zeta\caG_0$ of an extended frame
$F_\zeta\in\hat\caE$.

\subsection{The Symes formula.} A standard argument (cf.\ \cite[\S 2]{BurP2}) using the factorisation
theorem \ref{thm:factor} shows that there is a
well-defined map, the \emph{Symes map}, 
\begin{equation}\label{eq:Symes}
\Phi:\Ltg\to\hat\caE;\quad \Phi(\xi)=\exp(z\xi)_E.
\end{equation}
The group $\LIG$ acts adjointly on $\Ltg$ and the following argument
can be imported wholesale from \cite{BurP1}.
\begin{prop}\cite[\S 2]{BurP1}\label{prop:intertwine}
For $g\in\LIG$ and $\xi\in\Ltg$
\begin{equation}
[\Phi(\Ad g\cdot\xi)] = [g\sharp\Phi(\xi)].
\end{equation}
Therefore the map
\[
\bar\Phi:\Ltg\to\hat\caE/\caG_0;\quad \xi\mapsto [\Phi(\xi)],
\]
intertwines the adjoint action of $\LIG$ with its dressing action on $\hat\caL$.
\end{prop}
The main result concerning the Symes map is the following one, which tells us that a based 
extended frame
which admits a polynomial Killing field can be constructed via the Symes map \eqref{eq:Symes}.
\begin{thm}\label{thm:frame}
Let $F_\zeta\in\hat\caE$ admit an adapted polynomial Killing field $\xi(z) =
\zeta^{-4d-2}\xi_{-4d-2}+\zeta^{-4d-1}\xi_{-4d-1}+\ldots$. 
Set $\eta = \zeta^{4d}\xi(0)\in\Ltg$. Then $[F_\zeta] = [\Phi(\eta)]$.
\end{thm}
To prove this we must import a result regarding Lax equations of the form
\begin{equation}\label{eq:Lax2}
d\xi = [\xi,(\zeta^{4d}\xi dz)_E],
\end{equation}
which define a pair of commuting flows on the finite dimensional vector space
\[
\{\xi\in\Lcg:\zeta^{4d+2}\xi\in\Lig\}\subset\Lhg.
\]
The fact that \eqref{eq:Lax2} does indeed define a pair of commuting flows follows from \cite{BurP1}
\textit{mutatis mutandis}. We require another result from \cite{BurP1} which follows equally easily
by trivial modification.
\begin{lem}\cite[Thm 2.5]{BurP1}
Suppose $\xi(z)$ is a solution to the Lax equation \eqref{eq:Lax2}. Define
$\alpha_\zeta\in\Omega^1_\C\otimes\Lhg$ by $\alpha_\zeta = (\zeta^{4d}\xi dz)_E$. Then $\alpha_\zeta$
satisfies the Maurer-Cartan equations \eqref{eq:MC}.
\end{lem}
\begin{proof}[Proof of theorem \ref{thm:frame}]
Since $\Phi(\eta) = \exp(z\eta)_E$ there exists $\chi:\C\to\LIB$ for which $\exp(z\eta) =
\Phi(\eta)\chi$. Left invariant differentiation of this gives
\[
\eta\, dz  = d\Phi(\eta).\Phi(\eta)^{-1} + \Ad\Phi(\eta)\cdot(d\chi.\chi^{-1}).
\]
This shows that
\[
\Phi(\eta)^{-1}.d\Phi(\eta) = (\Ad\Phi(\eta)^{-1}\cdot\eta dz)_E =
(\Ad\Phi(\eta)^{-1}\cdot\zeta^{4d}\xi(0) dz)_E.
\]
Therefore $\tilde\xi = \Ad\Phi(\eta)^{-1}\cdot\xi(0)$ 
satisfies the Lax equation \eqref{eq:Lax2} and has initial condition $\tilde\xi(0)=\xi(0)$. By
uniqueness of solutions to o.d.e's, $\tilde\xi = \xi$ and therefore
$\Phi(\eta)^{-1}.d\Phi(\eta)=\alpha_\zeta$ with $\Phi(\eta)|_{z=0}=I$. Hence
$[\Phi(\eta)] = [F_\zeta]$.
\end{proof}

\subsection{Dressing orbits of the vacuum solutions.} To prove theorem \ref{thm:dressing} we
follow closely the lemmata presented in \cite{BurP2}, adapting them as we need. The first lemma
combines the appropriate versions of lemma 3.1 and proposition 3.2 from \cite{BurP2}. The proofs are
essentially the same and rely only on the properties that $A\in\fg$ be semisimple and an
eigenvector of $\tau$.
\begin{lem}
Let $\xi,\eta\in\Ltg$. Then $[\Phi(\xi)]=[\Phi(\eta)]$ if and only if
\begin{equation}\label{eq:etaxi}
\xi_{-2} = \eta_{-2},\quad \xi_{-1}=\eta_{-1},\quad \ad^n\eta\cdot\xi\in\Lig,\ \forall n\geq 1.
\end{equation}
Hence $[\Phi(\xi)] = [F^\mu]$ if and only if $\xi_{-2}=A$, $\xi_{-1}=0$ and $[\xi,A]=0$.
\end{lem}
The next result is the modification of \cite[Prop.\ 3.3]{BurP2} appropriate for our purposes.
\begin{lem}\label{lem:xieta}
For $\xi\in\Ltg$ and $g\in\LIB$ we have $[\Phi(\xi)]=g\sharp[F^\mu]$ if and only if $\Ad
g(0)^{-1}\cdot\xi_{-2}=A$ and $[\Ad g^{-1}\cdot\xi,A]=0$.
\end{lem}
\begin{proof}
By proposition \ref{prop:intertwine} we require conditions for which $[\Phi(\Ad g^{-1}\cdot\xi)] =
[F^\mu]$. By the previous lemma we obtain this precisely when
\[
(\Ad g^{-1}\cdot\xi)_{-2} = A,\quad (\Ad g^{-1}\cdot\xi)_{-1}=0,\quad  [\Ad g^{-1}\cdot\xi,A]=0.
\]
We can expand $\Ad g^{-1}\cdot\xi$ as a Laurent series in $\zeta$, using $g=(I+g_1\zeta +\ldots)g_0$
to obtain
\begin{equation}\label{eq:terms}
\Ad g^{-1}\cdot\xi = \Ad g_0^{-1}\left(\xi_{-2}\zeta^{-2} + ([\xi_{-2},g_1]+\xi_{-1})\zeta^{-1} + 
\ldots\right).
\end{equation}
The conditions we require now become
\begin{equation}\label{eq:cond}
\Ad g_0^{-1}\cdot\xi_{-2}=A,\quad [\xi_{-2},g_1]+\xi_{-1}=0,\quad  [\Ad g^{-1}\cdot\xi,A]=0.
\end{equation}
But now we observe that the last of these equations implies $[(\Ad g^{-1}\cdot\xi)_{-1},A]=0$, and since 
$\ker(\ad A)\cap \fg_{-1}=\{0\}$, this implies $(\Ad g^{-1}\cdot\xi)_{-1}=0$. Now an examination of
the coefficient of $\zeta^{-1}$ in \eqref{eq:terms} shows that the middle equation in
\eqref{eq:cond} is a consequence of the last.
\end{proof}
We are now in exactly the situation to which the arguments of
\cite[p372-373]{BurP2} apply, after suitable elementary modifications to account for the quadratic
dependence upon $\zeta$, to yield the following result. 
\begin{prop}\label{prop:dressing}
Let $\alpha_\zeta$ be an extended Maurer-Cartan form for which \eqref{eq:Maslov} holds
and let $\xi_\zeta$ be an adapted polynomial Killing field
of degree $4d+2$ for $\alpha_\zeta$. Set $\eta(\zeta) = \zeta^{-4d}\xi_\zeta(0)$. 
Then there exists $\e >0$ and $g\in\LIB$ for which $g\sharp[F^\mu]=[\Phi(\eta)]$.
\end{prop}
Note that to prove this, by \eqref{eq:Maslov} and the previous proposition, it suffices to construct 
$g\in\LIB$ for which
\begin{equation}\label{eq:gcond}
\Ad g(0)^{-1}\cdot A = A,\quad [A,\Ad g^{-1}\eta]=0.
\end{equation}
We leave it to the reader to check that the essential parts of \cite[p372-373]{BurP2} can be easily
adapted for this purpose. 

Now we may deduce theorem \ref{thm:dressing} by the following chain of reasoning. 
Given a HSL torus $f$ we
obtain an extended Maurer-Cartan frame $\alpha_\zeta$ which possesses a \pk field $\xi_\zeta$ and
therefore is has an extended Lagrangian frame of the form $\Phi(\eta)$. From the previous proposition
we obtain $g\in\LIB$ and a gauge transformation $k\in\caG_0$ for which $g\sharp F^\mu = 
\Phi(\eta)k$. Set $F_\zeta=\Phi(\eta)k$, then this is an extended Lagrangian frame for $f$ for which
$(gF^\mu)_E=F_\zeta$. Equation \eqref{eq:chi} follows. 
\begin{rem}\label{rem:orbits}
(i) Denote by $\caO^\e(\mu)\subset\hat\caL(\mu)$ the $\LIB$-orbit of the vacuum solution
$[F^\mu]$. It follows from proposition \ref{prop:intertwine} that $\caO^\e(\mu)\simeq 
\LIB/Z^\e_A$ where 
\begin{equation}\label{eq:Gamma}
Z^\e_A = \{ g\in\LIB: \Ad g\cdot \zeta^{-2}A = \zeta^{-2}A\}.
\end{equation}
Notice that this is independent of $\mu$.\\
(ii) Except for the proof of lemma \ref{lem:xieta}, everything we have done applies to the more 
general case of
an extended Maurer-Cartan form of degree $2$ which is $\tau$-equivariant for some finite order 
automorphism $\tau$ of order $k\geq 4$ on a reductive Lie algebra $\fg^\C$. In proposition
\ref{prop:dressing} we should explicitly require $\eta_{-2}$
to be semisimple, then the role of $A\in\fg_{-2}$ is played by the unique element in the $\Ad
B$-orbit of $\eta_{-2}$ for which $[A,\bar A]=0$: this provides a corresponding vacuum
solution. To prove lemma \ref{lem:xieta} we require the extra condition that
$\ker(\ad A)\cap \fg_{-1}=\{0\}$, which is true for our example but will not hold universally. In
fact this is the only obstruction in extending proposition \ref{prop:dressing} (and therefore
\cite[Thm 3.7]{BurP2}) to the more general setting.
\end{rem}

\end{document}